\documentclass[12pt,reqno]{amsart} 

\usepackage{etoolbox} 
\patchcmd{\subsubsection}{\normalfont}{\bfseries\textnormal}{}{} 

\usepackage[vmargin=2cm,hmargin=2cm]{geometry} 

\usepackage{caption} 
\usepackage{mathtools}

\usepackage{url}
\usepackage[english]{babel}
\usepackage{enumerate} 
 
\usepackage{bm} 
\usepackage[full]{textcomp}
\usepackage{newtxtext} 
\usepackage[cal=boondoxo]{mathalfa} 
\usepackage[bigdelims,vvarbb]{newtxmath} 

\usepackage{microtype} %

\newcommand{\dx}{\mathrm{d}}
\renewcommand{\Re}{\textit{Re}}

\newcommand{\scrA}{\mathscr{A}} 
\newcommand{\scrB}{\mathscr{B}} 
\newcommand{\scrC}{\mathscr{C}} 
\newcommand{\scrE}{\mathscr{E}} 
\newcommand{\scrO}{\mathscr{O}}

\newcommand{\C}{\mathbb{C}}

\newcommand{\N}{\mathbb{N}} 
\newcommand{\R}{\mathbb{R}}

\newcommand{\acc}{\mathbb{a}}

\newcommand{\scrF}{\mathcal{F}}
\newcommand{\boldF}{\bm{F}}

\newcommand{\Odi}[1]{\Odipm{}{#1}}

\newcommand{\Odipm}[2]{\boldsymbol{O}_{#1} (#2)}

\renewcommand{\qedsymbol}{$\square$}
\newenvironment{Proof}[1][Proof]{\par\noindent\textbf{#1.}~}
{\hfill\qedsymbol\smallskip\par}

\newtheoremstyle{itTheorems}
{10pt}
{6pt}
{\itshape}
{}
{\bfseries}
{.}
{.5em}
{\thmname{#1}\thmnumber{ #2}\thmnote{ (#3)}}

\theoremstyle{itTheorems} 
\newtheorem{Theorem}{Theorem}
\newtheorem{Proposition}{Proposition}
\newtheorem{Lemma}{Lemma}

\newtheorem{Definition}{Definition}
\newtheorem{Example}{Example}

\newtheoremstyle{remark}
{10pt}
{6pt}
{\rm} 
{}
{\bfseries}
{.}
{.5em}
{\thmname{#1}\thmnumber{ #2}\thmnote{ (#3)}}
\theoremstyle{remark} 
\newtheorem{Remark}{Remark}

\allowdisplaybreaks

\begin{document} 

\title[A unified strategy to compute some special functions]{A unified strategy to compute some special functions\\
of number-theoretic interest} 

\author{Alessandro Languasco}  

\begin{abstract}    
We introduce an algorithm to compute the functions belonging to a 
suitable set  $\scrF$ defined as follows: $f\in \scrF$   
means that $f(s,x)$, $s\in A\subset \R$ being fixed and  $x>0$, has a power series  
expansion centred at $x_0=1$ with convergence radius greater or equal than $1$; moreover,
it satisfies a functional equation of step $1$
and the Euler-Maclaurin summation formula can be applied to $f$.
Denoting the Euler gamma-function as $\Gamma$,
we will show that, for $x>0$,  $\log \Gamma(x)$,  the digamma
function $\psi(x)$, the polygamma functions $\psi^{(w)}(x)$, $w\in \N$, $w\ge1$, and, 
for $s>1$ being fixed, the Hurwitz $\zeta(s,x)$-function and its first partial derivative
$\frac{\partial\zeta}{\partial s}(s,x)$ are in $\scrF$.
In all these cases the coefficients of the involved power series will depend
on the values of $\zeta(u)$, $u>1$, where $\zeta$ is the Riemann zeta-function.
As a by-product, we will also show how to compute the Dirichlet $L$-functions
$L(s,\chi)$ and $L^\prime(s,\chi)$, $s> 1$, $\chi$ being a primitive Dirichlet character,
by inserting the reflection formulae of  $\zeta(s,x)$ and $\frac{\partial\zeta}{\partial s}(s,x)$
into the first step of the Fast Fourier 
Transform algorithm. Moreover, we will obtain some new formulae and algorithms for the 
Dirichlet $\beta$-function and for the Catalan constant $G$. 
Finally, we will study the case of the Bateman $G$-function
and of the alternating Hurwitz zeta-function, also known as the $\eta$-function;
we will show that, even if they are not in $\scrF$, our approach
can be adapted to handle them too.
In the last section we will also describe some tests 
that show a performance gain with respect to a standard multiprecision
implementation of  $\zeta(s,x)$ and $\frac{\partial\zeta}{\partial s}(s,x)$, $s>1$, $x>0$.
\end{abstract}  

\maketitle 
\makeatletter
\def\subsubsection{\@startsection{subsubsection}{3}%
\z@{.3\linespacing\@plus.5\linespacing}{-.5em}%
{\normalfont\bfseries}} 
\makeatother

\section{Introduction}
The goal of this paper is to show that for a suitable
set of functions $\scrF$  there exists a unified computational
strategy, at least when the main variable is a positive
real number.  
The set of functions we will work on is described in the following 
\begin{Definition}[The set $\scrF$]
\label{set-def}
We will say that a function $f\colon A\times (0,+\infty) \to \R$, where 
$A\subset \R$, is in the set $\scrF$ if and only if 
it has the following properties:
\begin{enumerate}[i)]
\item
\label{series-f}
for every fixed $s\in A$
there exists a sequence of functions $c_{f}(s,k)\in \R$ and $\rho_f(s)\ge 1$
such that $\sum_{k=0}^{\infty} c_{f}(s,k)(1-z)^k$ is absolutely convergent
to  $f(s,z)$ for every  $z\in(1-\rho_f(s), 1+\rho_f(s))$;
\item 
\label{control-number-summands-f}
for every fixed $s\in A$,
there exist $k_{f}(s) \in \N$ and $C_{f}(s,k)> 0$ such that 
$\vert c_{f}(s,k) \vert \le C_f(s,k)$ and $C_f(s,k)$ is a decreasing sequence for every $k\ge k_{f}(s)$;
\item 
\label{recursive-f}
for every fixed $s\in A$, $f(s,\cdot)$ satisfies a functional equation of step $1$, i.e, 
there exists 
a function $g_f\colon A\times (0,+\infty) \to \R$ such that, for every fixed $s\in A$, we have  
\begin{equation}
\label{difference-f}
f(s,z+1)=f(s,z)+ g_f(s,z) \quad \text{for every}\ z\in (0,+\infty);
\end{equation}
\item 
for every fixed $s\in A$,
the function $g_f(s,\cdot)$ in \eqref{difference-f} can be used in the Euler-Maclaurin 
summation formula; i.e, for  every fixed $s\in A$, $g_f(s,\cdot)$ verifies the 
hypotheses of Lemma \ref{Euler-Maclaurin} below.
\end{enumerate} 
\end{Definition}

As we will see later, $\log\Gamma(x) \in \scrF$ for $x>0$, ($\Gamma$ is the \emph{Euler gamma-function})
together with the \emph{digamma function} $\psi(x)$, and the treatment here presented 
generalizes and improves on the one in \cite{Languasco2021}. These functions have the 
coefficients $c_{f}(s,k)$ depending on the values of the \emph{Riemann zeta-function}
at positive integers greater than $1$; thus, to evaluate them for $x\in (0,1)$, 
we will need to precompute a sufficiently large number of values of $\zeta(j)$, $j\in \N$, $j \ge 2$. 
For $j$ even, 
the Bernoulli numbers are hence needed.
For $x>1$ we will use the Euler-Maclaurin formula, and hence we will need  
to precompute a sufficiently large set of Bernoulli numbers; to be able to obtain
a sufficiently good accuracy in this step, we will also need to exploit an
\emph{``horizontal shift''} trick, 
see Section \ref{horiz-shift}, which will require to compute some values of the 
function $g_f(s,x)$ used in point \ref{recursive-f}) of Definition \ref{set-def}.  
If $x$ is large and some asymptotic formula for $f(s,x)$, $s$
being fixed, $x\to +\infty$, is known, its use might require a smaller computational effort 
than the one needed to perform the recursive step of point \ref{recursive-f})
of Definition \ref{set-def}.

Another special function belonging to $\scrF$ is
the \emph{Hurwitz zeta-function} $\zeta(s,x)$,
together with its first partial derivative
$\zeta^\prime(s,x):=\frac{\partial\zeta}{\partial s}(s,x)$,
 for $s>1$ being fixed and  $x>0$.
In these cases the coefficients $c_{\zeta_{H}}(s,k)$, $c_{\zeta^\prime_{H}}(s,k)$ will depend on
$\zeta(w),\zeta^{\prime}(w)$, $w>1$, and on some \emph{Euler beta-function} values. 
All the  functions mentioned before will be involved in 
computing $L(s,\chi)$ and $L^{\prime}(s,\chi)$, $s\ge1$, 
where $L$ denotes a \emph{Dirichlet $L$-function} and
$\chi$ is a non-principal  Dirichlet character modulo an odd prime number $q$. 

As a by-product, since the \emph{polygamma} functions $\psi^{(w)}(x)$, $w\in \N$, $w\ge 1$,
can be written in terms of $\zeta(w+1,x)$, we will also obtain
that $\psi^{(w)}(x)\in \scrF$, for every $x>0$ and $w\in \N$, $w\ge 1$. 
Moreover, we will show some new formulae and algorithms for the 
\emph{Dirichlet $\beta$-function} $\beta(s)$, together with $\beta^\prime(s)$ 
and $\beta^\prime(s)/\beta(s)$, $s>1$, too; we will also obtain a new fast convergent 
series for the \emph{Catalan} constant $G$. 
Finally, we will  discuss  two other examples regarding the 
\emph{Bateman} $G$-function and  the \emph{alternating Hurwitz zeta-function}, 
also known as the $\eta$-function,
to show that, even if $G(\cdot),\eta(s,\cdot)\not \in \scrF$ because
they are not the solution of a functional equation of the type in 
point  \ref{recursive-f}) of Definition \ref{set-def}, 
our strategy can be adapted to these cases too.

In the final section we will report on some practical experiments.  
From all these tests and examples, we can infer that the algorithm here presented
is particularly useful when it is possible to exploit the precomputation
of the $c_f(s,k)$-coefficients of the series in point \ref{series-f}) of Definition \ref{set-def}, 
as for the Dirichlet $L$-functions $L(s,\chi)$, $s\ge1$, and its
first derivative, where $\chi$ is a non-principal Dirichlet character mod $q$, and $q$ runs 
into a large set of odd primes. 

\subsection*{Outline.}
In Section \ref{lemmas} we will collect some useful lemmas, and in the following
one we will write the algorithm for $f\in \scrF$; we will also study its computational
cost. 
Section  \ref{gamma-section} is dedicated to show a detailed treatment for the $\log\Gamma$-function,
while in Section \ref{digamma-section} we will discuss the case of the digamma function $\psi(x)$.
Section \ref{Hurwitz-zeta-section} is devoted to show how to adapt the general algorithm
to the case of the Hurwitz zeta-function and, as a by-product, to the polygamma functions.
Section \ref{Hurwitz-prime-section} is about $\zeta^\prime(s,x)$,
$x>0$ and $s>1$, since it will be useful to handle the first derivative of the
Dirichlet $L$-functions. In the same section, we will also treat the case of the 
Dirichlet $\beta$-function and of  the Catalan constant $G$.
In Section \ref{DirichletL}
we will then discuss how to use the reflection formulae of
the previously mentioned functions in the Fast Fourier Transform algorithm with the goal of computing
the Dirichlet $L$-functions attached to non-principal  primitive Dirichlet characters.
In Section \ref{beyond-F} we will describe how to handle the cases of
the Bateman $G$-function and the alternating Hurwitz zeta-function even
if they are not in $\scrF$.
Moreover,  to show the good performances of this algorithm, 
in Section \ref{tests} we will describe some results
obtained with programs developed using 
Pari/GP \cite{PARI2021}; source codes and examples are available
on the following page: \url{http://www.math.unipd.it/~languasc/specialfunctions.html}.
%

\subsection*{Acknowledgment.}
I would like to thank the anonymous referee for his/her remarks
and suggestions.
I would also thank my colleague and friend Mauro Migliardi
for having read a preliminary version of this paper.

\section{Lemmas}
\label{lemmas}
To estimate the error we have in approximating $f$ with a finite sum
in point  \ref{series-f}) of Definition \ref{set-def},
we will need the following lemma
which is an adaptation of the ratio test.
\begin{Lemma} 
\label{stima-errore-criterio-rapporto}
Let $c_k\in \C$,  $\sum_{k=0}^{\infty}c_k$  be an absolutely convergent series,
$\mu\in(0,1)$ and $K\in \N$ such that
\(
\vert
c_{k+1}/c_k
\vert 
\le \mu
\) for every $k \ge K$. 
Denoting $\sum_{k=0}^{\infty}c_k =: S \in \C$ and $\sum_{k=0}^{m}c_k =: S_m \in \C$, 
we have that
\(
\vert S - S_m \vert \le \vert c_{m+1} \vert/(1-\mu)
\)
for every \( m \ge K-1\).
\end{Lemma} 
\begin{Proof}
Let $m+1\ge K$. We have 
$\vert c_{m+2}\vert  \le \mu \vert c_{m+1}  \vert$  
and
$\vert c_{m+3}\vert  \le \mu \vert c_{m+2}\vert \le \mu^2\vert c_{m+1}\vert$.
Arguing by induction, we prove, for every $p\in \N$, $p\ge 1$, that
\(
\vert c_{m+p}\vert  \le \mu^{p-1}\vert c_{m+1}\vert.
\)
Hence we get
\[
\vert S - S_m \vert
=
\Bigl\vert \sum_{k=m+1}^{\infty} c_k  \Bigr\vert
\le
\sum_{k=m+1}^{\infty} \vert c_k  \vert
= 
\sum_{p=1}^{\infty} \vert c_{m+p}  \vert
\le
\sum_{p=1}^{\infty} \mu^{p-1}\vert c_{m+1}  \vert
\le
\frac{\vert c_{m+1}  \vert}{1-\mu} 
\]
using the well known theorem about the geometric series.
\end{Proof}

A similar result holds true for the root test too.
We will also need a statement about the Euler-Maclaurin formula; we will use the following one
that can be obtained by combining the topics in Sections 3.3-3.4  of 
Stoer-Bulirsch \cite{StoerBulirsch2002}, or referring to 
Cohen \cite[Corollary 9.2.3(2) and Proposition 9.2.5(2)]{Cohen2007}. 
\begin{Lemma}
\label{Euler-Maclaurin}
Let $a\in \R$, $N\in \N$, and $m\in \N$, $m\ge 1$.
Assume that $h\in C^{2m+4}([a, a+ N])$ and that
both $h^{(2m+2)}, h^{(2m+4)}$ have constant sign on $[a,a+N]$.
We have
\begin{align*}
\sum_{j=0}^{N} h(a+j)  
&
= 
\int_a^{a+ N} h(w)\ \dx w
+ 
\frac12  \bigl( h(a+ N) + h(a)  \bigr)
+
\sum_{n=1}^{m}  \frac{B_{2n}}{(2n)!} \bigl(h^{(2n-1)}(a+ N)- h^{(2n-1)}(a) \bigr)
\\&
- \frac{1}{(2m)!}  \int_a^{a+ N}  B_{2m}(\{w-a\}) h^{(2m)}(w)\ \dx w,
\end{align*}
where $B_n(u) $ are the Bernoulli polynomials 
and $B_{n}$ are the  Bernoulli numbers of order $n$.
Moreover\
\[
\Bigl \vert
\frac{1}{(2m)!}  \int_a^{a+ N}   B_{2m}(\{w-a\}) h^{(2m)}(w)\ \dx w
\Bigr \vert
\le 
\frac{\vert B_{2m+2} \vert}{(2m+2)!} \Bigl \vert h^{(2m+1)}(a+ N)-h^{(2m+1)}(a) \Bigr\vert.
\]
\end{Lemma}

We will also need the following elementary estimates about the digamma function and 
the Riemann zeta-function for $x>0$.
\begin{Lemma}
\label{elementary-estim}
Let   $\psi(\cdot)$ be the digamma function, $\zeta(\cdot)$ be the Riemann zeta-function and 
let $x>0$. Then 
\(
\log x - 1/x < \psi(x) < \log x.
\)
Moreover, for every $x>1$, we have that
\[
1+ \frac{1}{2^{x}}<\zeta(x) < 1+\frac{1}{2^{x}}\frac{x+1}{x-1}
\]
and, for $x\ge 3$, also that
\[
-\frac{\log2 +(2/3)\log 3}{2^x} < \zeta^{\prime}(x) <-\frac{\log 2}{2^x}.
\]
\end{Lemma}
\begin{Proof}
The first inequality follows from Theorem 5 of Gordon \cite{Gordon1994}.
The estimate on $\zeta(x)$, $x>1$, can be obtained  from the definition
of the Riemann zeta-function in  $\Re(u) >1$, and the integral test for the series.
Recalling that
$\zeta^{\prime}(u) = -\sum_{n=2}^{\infty} (\log n)  n^{-u}$,  $\Re(u) >1$,
we have $-\zeta^{\prime}(x) > (\log 2) 2^{-x}$ for every $x> 1$. Moreover, using
that $(\log x)/x$ is a decreasing sequence for  $x\ge e$,
the last part of the lemma  follows by  remarking
\[
-\zeta^{\prime}(x) =  \frac{\log 2}{2^x} + \sum_{n=3}^{\infty} \frac{\log n}{ n^{x}}
<
\frac{\log 2}{2^x} +  \frac{\log 3}{3} \sum_{n=3}^{\infty} \frac{1}{n^{x-1}}
= 
\frac{\log 2}{2^x} +  \frac{\log 3}{3} \Bigl(\zeta(x-1)-1- \frac{1}{2^{x-1}}\Bigr)
\]
and  using the inequality $\zeta(x) <1+2^{1-w}$, $x\ge 3$, which 
follows from the one previously proved.
\end{Proof}

\section{The Algorithm for $f\in \scrF$}
Let $f\in \scrF$, $s\in A$ be fixed and  $x>0$. Denote by $\lfloor x \rfloor$ and $\{ x\}$
the integral and fractional parts of $x$.
If $\lfloor x \rfloor   \ge 1$ and  $\{ x\} >0$, 
we can use \eqref{difference-f} 
to write
\begin{equation} 
\label{main-recursive-f}
f(s,x)
=
f(s, x -1)
+ 
g_f(s, x -1)
= 
\dotsm
=
f(s, \{ x\})
+
\sum_{j=0}^{ \lfloor x\rfloor -1} 
g_f(s, \{ x\} + j).
\end{equation}
In the case in which  $x=\ell \in \N$, $\ell\ge 2$,
formula \eqref{main-recursive-f} becomes
\begin{equation}
\label{main-recursive-f-simplified}
f(s,\ell)
=
f(s, \ell -1) + g_f(s, \ell -1)
= 
\dotsm
=
f(s,1) 
+
\sum_{j=1}^{ \ell -1} 
g_f(s,j) .
\end{equation}
In this way we split the solution in two parts: the first evaluates $f(s, \{ x\})$ 
and the second computes the \emph{tail}, i.e,
the difference $f(s,x)-f(s, \{ x\})$
expressed as a  sum of values of the $g_{f}$-function. Unfortunately the length
of such a sum may be too large, so we will evaluate
the tail using Lemma \ref{Euler-Maclaurin}.  
We see now how to perform the computation of the quantities in \eqref{main-recursive-f};
the ones in \eqref{main-recursive-f-simplified} can be handled analogously.

\subsection{Computation of $f(s, \{ x\})$, $\{ x\} > 0$, $s\in A$ being fixed} 
\label{comput-f-in-(0,1)}
Recalling that  $s\in A$ is fixed, 
point \ref{series-f}) of Definition \ref{set-def}
gives 
\begin{equation}
\label{f-series}
f(s,z) =\sum_{k=0}^{\infty} c_{f}(s,k)(1-z)^k,
\end{equation}
and the series in \eqref{f-series}
absolutely converges for $z\in(1-\rho_f(s),1+\rho_f(s))$.
Since $ \rho_f(s) \ge 1$,
this fact and \eqref{difference-f} let us obtain $f(s, \{ x\})$ directly
or as $ f(s, \{x\}+1) - g_f(s, \{ x\})$.
In this way  we will be able to have a good estimate
on the number of  summands needed to approximate, up to the desired accuracy, 
the  series in \eqref{f-series} with a finite sum.

Before being more precise, we need the following remark.
Using point \ref{control-number-summands-f})
of Definition \ref{set-def}, we have 
$\frac{C_{f}(s,k+1)}{C_{f}(s,k)}  \le 1$ for $k \ge k_{f}(s)$
and hence,
using Lemma \ref{stima-errore-criterio-rapporto}
on  \eqref{f-series}, we have 
$\mu_f = \vert 1-z\vert$. 
Letting now $\acc\in \N$, $\acc\ge 2$, we have that
for every $z\in  (0,2)\subset(1-\rho_f(s),1+\rho_f(s))$ 
there exists $r=r_f(s,z, \acc ) \ge k_{f}(s)$
such that 
\begin{align}
\notag
\Bigl\vert \sum_{k=r+1}^{\infty} c_{f}(s,k)(1-z)^k \Bigr\vert
&\le
\sum_{k=r+1}^{\infty} C_{f}(s,k)\vert 1-z\vert^k
\le 
C_{f}(s,r+1) 
\frac{\vert 1-z \vert^{r+1}}{1- \vert 1-z\vert} 
\\&
\label{tail-f-series}
\le 
C_{f}(s, k_{f}(s)+1) 
\frac{\vert 1-z \vert^{r+1}}{1- \vert 1-z\vert} 
<
2^{-\acc-1},
\end{align}
since $C_{f}(s,k)> 0$ is a decreasing sequence for $k\ge k_{f}(s)$.
A straightforward computation reveals that in \eqref{tail-f-series} we can choose 
\begin{equation}
\label{r-gen-estim}
r_f(s,z,\acc) = \max
\Bigl(
\Bigl\lceil \frac{(\acc+1) \log 2 + \vert \log (1-\vert 1-z\vert) \vert 
+ \vert \log ( C_{f}(s,k_{f}(s)+1) ) \vert}
{ \vert \log  \vert 1-z \vert \vert}\Bigr\rceil -1;
k_{f}(s)
\Bigr),
\end{equation}
where we denoted  as $\lceil u \rceil$ 
the least integer greater than or equal to 
$u\in \R$.

\subsubsection{Optimising the number of summands for  $f(s,\{x\})$, $s\in A$ being fixed}
\label{shifting-trick-f}
Clearly $r_f(s,z,\acc) $ becomes larger as $\vert 1-z\vert$ increases. So
when $0<z<1/2$ we will evaluate $f$ at $1+z$
exploiting \eqref{difference-f}. 
In this way we will always use the best convergence interval, $z\in(1/2,3/2)$, we have for
the series in \eqref{f-series}; we also remark that 
\[
r_f(s,z,\acc) \le r_f(s,1/2,\acc) =r_f(s,3/2,\acc) =  \max\Bigl(\acc+1 +
\Bigl\lceil
\frac{\vert\log ( C_{f}(s,k_{f}(s)+1) ) \vert}{\log 2}
\Bigr\rceil;
k_{f}(s)
\Bigr)
\]
for every $z\in(1/2,3/2)$.  
Summarising,  using  \eqref{f-series}-\eqref{tail-f-series},
for $\{x\}\in(1/2,1)$ we have that there exists $\theta_1=\theta_1(s,\{x\}) \in (-1/2,1/2)$ 
such that
\begin{equation}
\label{f-x>1/2}
f(s, \{ x\})
=
\sum_{k=0}^{\infty} c_{f}(s,k)(1-\{x\})^k
= 
\sum_{k=0}^{r_f(s,\{x\},\acc)}  c_{f}(s,k)(1-\{x\})^k
+\vert \theta_1 \vert 2^{-\acc}.
\end{equation}
We also  remark that for  $\{x\}\in(1/2,1)$, we have 
\begin{equation} 
r_f(s,\{x\},\acc)  
\label{r-f-x<1/2} 
\le 
\max
\Bigl(\acc+1 + 
\Bigl\lceil\frac{\vert \log ( C_{f}(s,k_{f}(s)+1) ) \vert}{\log 2}\Bigl\rceil;
k_{f}(s)
\Bigr).
\end{equation}
Moreover, using \eqref{difference-f} and \eqref{f-series}-\eqref{tail-f-series}, 
for $\{x\}\in(0,1/2)$ we have  that there exists 
$\theta_2=\theta_2(s,\{x\}) \in (-1/2,1/2)$ such that
\begin{align}
\notag
f(s, \{ x\}) &=  - g_f(s, \{ x\}) + f(s,\{x\}+1)
=
- g_f(s, \{ x\})   + \sum_{k=0}^{\infty}  c_{f}(s,k) (-\{x\}) ^k
\\
\label{f-x<1/2}
&=
- g_f(s, \{ x\})  +  \sum_{k=0}^{r^\prime_f(s,\{x\},\acc)} c_{f}(s,k) (-\{x\}) ^k
+\vert \theta_2 \vert 2^{-\acc},  
\end{align}
where 
\begin{align}
\notag
r^\prime_f(s,\{x\},\acc)&:=r_f(s,1+\{x\},\acc) 
\\&\label{rprime-f-x<1/2}
=  \max
\Bigl(
\Bigl\lceil \frac{(\acc+1) \log 2 + \vert \log ( 1- \{x\} )\vert 
+ \vert\log ( C_{f}(s,k_{f}(s)+1) )\vert}{ \vert \log \{x\} \vert} 
\Bigr\rceil-1;
k_{f}(s)
\Bigr)
\end{align} 
and
\[
r^\prime_f(s,\{x\},\acc)
\le
\max
\Bigl(\acc+1 + \Bigl\lceil\frac{ \vert\log ( C_{f}(s,k_{f}(s)+1) )\vert}
{ \log 2}\Bigl\rceil;
k_{f}(s)
\Bigr).
\]
Assuming that the needed $c_{f}(s,k)$-values  can be precomputed and stored 
with a sufficiently good accuracy,
the formulae in \eqref{f-x>1/2} and \eqref{f-x<1/2} allow us to compute
$f(s,\{x\})$, with an accuracy of $\acc$ binary digits 
using about $\acc+1$ summands.
We will be more precise on this in the next subsection. 

\begin{Remark} For $\lfloor x \rfloor  \ge 1$, $\{x\}\in(0,1/2)$, combining
the first part of \eqref{f-x<1/2} with \eqref{main-recursive-f}
we obtain
\begin{equation} 
\label{main-recursive-f-alt}
f(s,x)  
=
f(s, \{ x\}+1)
+
\sum_{j=1}^{ \lfloor x\rfloor -1} 
g_f(s, \{ x\} + j).
\end{equation}
It is in fact slightly faster to use \eqref{main-recursive-f-alt} since
it involves a shorter sum than the one in the right hand side of  \eqref{main-recursive-f}
and $\{ x\}+1 \in (1,3/2)$ is already in the best convergence
interval for the series that defines $f(s,\cdot)$.
\end{Remark}

\subsection{Computational costs and error terms for $f(s,x), x\in(0,1)$, $s\in A$ being fixed}
\label{comp-cost-f} 
We assume the  value of $f(s,1/2)$ is known.
The estimates on $r_f(s,\{x\},\acc),r^\prime_f(s,\{x\},\acc)$ in \eqref{r-f-x<1/2} 
and \eqref{rprime-f-x<1/2} for
every $\{x\}\in (1/2,1)$ and, respectively, $\{x\}\in (0,1/2)$, imply
that $f(s, \{ x\})$, $x\in (0,1)$, can be obtained with an $\acc$-bit
accuracy using about
$\acc+1$ summands.
The summation can be performed combining 
the \emph{pairwise summation} \cite{Higham1993}  algorithm
with Kahan's \cite{Kahan1965} method  (the minimal block
for the pairwise summation algorithm is summed using 
Kahan's procedure)
to have a good compromise between accuracy, computational
cost and execution speed.
Since the needed powers of $\{ x\}$ and $1-\{ x\}$ can be obtained using a repeated product strategy,
the cost of computing $f(s, \{ x\})$, $\{x\}\in (1/2,1)$
is essentially $\Odi{\acc}$ evaluations of $c_{f}(s,k)$,  
$\Odi{\acc}$ floating point products, $\Odi{\acc}$ floating point summations
with an accuracy of $\acc$ binary digits; 
for $\{x\}\in (0,1/2)$, the additional cost of computing $g_f(s , \{ x\})$
must be considered.

\subsection{The tail}
We explain here the general idea on how
to evaluate the tail of the problem in \eqref{main-recursive-f}, i.e,
\(
f(s,x)
-
f(s, \{ x\})
=
\sum_{j=0}^{ \lfloor x\rfloor -1} 
g_f(s, \{ x\} + j).
\)
Using  Lemma \ref{Euler-Maclaurin} with $h(\cdot) = g_f(s,\cdot)$, $s\in A$ being fixed, 
the error terms will depend on the size of
the derivatives of $h$.
Here we  just make some general considerations, but
the actual evaluation of the error term of this part will be directly performed
for the examples we will see in the next sections.
We just recall here that, in general, a well known characteristic of the 
Euler-Maclaurin formula, see, e.g., Stoer-Bulirsch \cite[Sections 3.3-3.4]{StoerBulirsch2002},
is that the error term might not converge to $0$; but,
in fact, it initially decreases until it reaches a minimal value. 
For an example of this phenomenon, see Example \ref{example-hor-shift} 
and Table \ref{EM-error-table}  of Section \ref{why-hor-shift}.
We will need then to evaluate if such a minimal value satisfies our desired accuracy; if this
is not the case, it is possible to improve on it by using the following idea.

\subsubsection{The horizontal shift}
\label{horiz-shift}
The tail-computation can be rearranged in the following way. Letting $t\in \N$,
$2\le t \le \lfloor x \rfloor$, and $v := \{ x\}+t> t$. Moreover, letting
$u:=\lfloor x \rfloor-t \in \N$, we can write that
\(
f(s, \lfloor x \rfloor + \{ x\})   
= 
f(s,u  +  v)   
\)
and hence
\begin{align}
\notag
f(s,\lfloor x \rfloor + \{ x\}) 
-
f(s,\{ x\})  
&=
f(s,u  +  v)  
-
f(s,v)  
+
f(s, v) 
-
f(s, \{ x\}) 
\\
\label{horiz-shift-f}
&=
f(s,u  +  v)    
-
f(s,v) 
+
\sum_{j=0}^{t-1} g_f(s,\{ x\}+ j ) .
\end{align}
Now we can apply  the Euler-Maclaurin formula to $f(s ,u  +  v)   - f(s, v)$
in a situation, $v>t \ge 2$, in which, as we will see in the next sections, we will 
have a small error term. 
Hence, by choosing $t$ large enough, and at the cost
of evaluating $t$-times the $g_f$-function, we can obtain a final better accuracy.

If $x$ is large and some asymptotic formula as $x \to +\infty$ is known for $f(s,x)$, $s$
being fixed, it might be better to replace the recursive step of point \ref{recursive-f})
of Definition \ref{set-def} with such an asymptotic formula.

\subsection{The reflection formulae for $f(s,x)$, $x\in(0,1)$, $s\in A$ being fixed}
In the applications involving the use of the Fast Fourier Transform, see Section \ref{FFT-setting},
we can speed up the  global computation, and improve on the memory usage,
by using the values of $f(s,x) \pm f(s, 1-x)$  for $x\in(0,1)$, 
$s\in A$ being fixed, instead of 
the ones of $f(s,x)$.
A nice feature of using  \eqref{f-x>1/2} and \eqref{f-x<1/2} for this goal
is that the odd or the even summands of the series of such quantities will vanish
thus reducing by a factor of $2$ the computational effort of getting such values.
We summarise the situation in the following
\begin{Proposition}
\label{DIF-formulae-f}
Let $f\in \scrF$, $s\in A$ being fixed, $x\in (0,1)$, $x\ne 1/2$,
$\acc\in \N$, $\acc\ge 2$.
Recalling \eqref{r-f-x<1/2} and \eqref{rprime-f-x<1/2}, let further 
$r_1(s,x,\acc) = r^\prime_f(s,x,\acc) /2$
and  $r_2(s,x,\acc) = r_f(s,x,\acc) /2$. 
There exists $\theta=\theta(s,x) \in (-1/2,1/2)$  such that
for $0<x  < 1/2$ we have
\begin{align} 
f(s,x)&+f(s, 1-x) 
\label{f-DIF-even-x<1/2}
= 
-  g_f(s,x) + 
2\sum_{\ell=0}^{r_1} c_{f}(s,2\ell) x^{2\ell}
+ \vert \theta \vert 2^{-\acc}, 
\\
f(s,x)&-f(s, 1-x) 
\label{f-DIF-odd-x<1/2}
= 
-  g_f(s,x) 
-2 
\sum_{\ell=0}^{r_1} c_{f}(s,2\ell+1) x^{2\ell+1}
+ \vert \theta \vert 2^{-\acc}, 
\end{align}
and  for $1/2<x<1$ we have
\begin{align}
f(s,x)&+f(s, 1-x) 
\label{f-DIF-even-x>1/2}
=
-  g_f(s, 1-x) + 2
\sum_{\ell=0}^{r_2} c_{f}(s,2\ell)  (1-x)^{2\ell}
+ \vert \theta \vert 2^{-\acc},
\\ 
f(s,x)&-f(s, 1-x) 
\label{f-DIF-odd-x>1/2}
=
g_f(s, 1-x) 
+2 \sum_{\ell=0}^{r_2}c_{f}(s,2\ell+1) (1-x)^{2\ell+1}
+ \vert \theta \vert 2^{-\acc}.
\end{align}
\end{Proposition}

\begin{Proof}
Assume that $0<x < 1/2$; in this case we compute $f(s,x)$ with \eqref{f-x<1/2}
and $f(s,1-x)$ with \eqref{f-x>1/2}. 
Since the involved series  absolutely converge, their sum is   the series of
$c_{f}(s,k)(x^k+(-x)^k) = 2 c_{f}(s,k) x^k$ when $k$ is even and zero otherwise.
Arguing as in Section \ref{shifting-trick-f} and remarking that  
$r_1(s,x,\acc)= r_f(s,1-x,\acc) /2=  r^\prime_f(s,x,\acc)/2$,
we  have that \eqref{f-DIF-even-x<1/2} holds.
Assume that $1/2<x<1$; in this case we compute 
$f(s,x)$ with \eqref{f-x>1/2}
and $f(s,1-x)$ with \eqref{f-x<1/2}. 
Arguing as before and remarking that  
$r_2(s,x,\acc)= r_f(s,x,\acc)/2 = r^\prime_f(s,1-x,\acc)/2$,
we  have that \eqref{f-DIF-even-x>1/2} holds.
In a similar way we can prove \eqref{f-DIF-odd-x<1/2} and  \eqref{f-DIF-odd-x>1/2};
the only difference is that now the summands are 
$c_{f}(s,k)(x^k-(-x)^k) = 2 c_{f}(s,k)x^k$ when $k$ is odd and zero otherwise.
This completes the proof.
\end{Proof}

We remark that in Proposition \ref{DIF-formulae-f} we have $r_2(s,x,\acc)
= r_1(s,1-x,\acc)$ for $x\in (0,1)$ and hence the right hand side of \eqref{f-DIF-even-x>1/2} 
can be obtained from the right hand side of \eqref{f-DIF-even-x<1/2} by replacing
any occurrence of $x$ with $1-x$ and vice versa.
Something similar hold for \eqref{f-DIF-odd-x<1/2} and \eqref{f-DIF-odd-x>1/2} if we also
perform a sign change.

We finally remark that the proof of Proposition \ref{DIF-formulae-f} reveals
that a similar statement holds for the infinite series 
of $f(s,x)\pm f(s,1-x)$  too.

\section{The $\log\Gamma$-function}
\label{gamma-section}

Our first result in this section is the following
\begin{Theorem}[The $\log \Gamma$-function is in $\scrF$]
$\log\Gamma(x)$, $x>0$, is in $\scrF$.
\end{Theorem}
\begin{Proof}
Recalling that  the well known Euler formula, see, e.g., Lagarias 
\cite[section 3]{Lagarias2013}, gives 
\begin{equation*}
\log\Gamma(z) = \gamma(1-z)  + \sum_{k=2}^{\infty} \frac{\zeta(k)}{k}(1-z)^k,
\end{equation*}
where $\zeta(\cdot)$ is the Riemann zeta-function, $\gamma$ 
is the Euler-Mascheroni constant and $z\in(0,2)$, by letting 
\begin{equation}
\label{coeff-gamma-def}
c_{\Gamma}(k) = \frac{\zeta(k)}{k}\ \textrm{for}\ k\ge 2,\quad c_{\Gamma}(1) = \gamma, \quad
c_{\Gamma}(0) =0 \quad 
\textrm{and}\quad C_{\Gamma}(k):=c_{\Gamma}(k)\ \textrm{for} \ k\ge 1,
\end{equation} 
we have that point \ref{series-f}) of  Definition \ref{set-def} 
holds. 
By Lemma \ref{elementary-estim} 
we have $1+ 2^{-w}<\zeta(w)<1+2^{1-w}$ for $w\ge 3$, so that
\begin{equation}
\label{coeff-gamma-decreasing}
\frac{c_{\Gamma}(k+1)}{c_{\Gamma}(k)} 
=
\frac{k \zeta(k+1)}{(k+1)\zeta(k)}
<
\frac{k}{(k+1)} \frac{1+2^{-k}}{1+2^{-k}}
<
1
\quad 
\text{for}
\ k\ge 3.
\end{equation}
Hence point \ref{control-number-summands-f}) of Definition \ref{set-def} also holds
with  $\overline{k}_\Gamma =3$. 
Using the well known relation $\Gamma(z+1) = z \Gamma(z)$, $z>0$, we obtain the functional equation  
\begin{equation}
\label{difference-gamma}
\log \Gamma(1+z) =\log \Gamma(z)  + \log z, \quad z>0,
\end{equation}
and hence we can say that point \ref{recursive-f}) of Definition \ref{set-def} holds 
with $g_{\Gamma}(s,z) = \log z$, $z>0$.
This proves that $\log\Gamma(x)$, $x>0$, is in $\scrF$.
\end{Proof}

The analysis on the number of summands that follows from Sections
\ref{comput-f-in-(0,1)}-\ref{comp-cost-f} gives the results already
proved in Section 3 of \cite{Languasco2021}; the same happens for the 
reflection formulae we obtain by specialising Proposition \ref{DIF-formulae-f}.
For completeness we insert here such results.
Recalling that   $\log \Gamma(1/2)=(\log \pi)/2$, we have the following
\begin{Proposition}[The number of summands for $\log \Gamma$]
Let $x\in (0,1)$, $x\ne 1/2$,
$\acc\in \N$, $\acc\ge 2$,
\begin{equation}
\label{r-gamma-estim}
r_{\Gamma}(x,\acc) = 
\max
\Bigl\{
\Bigl\lceil \frac{(\acc+1) \log 2 
+ \vert \log (1-\vert 1-x\vert)\vert}{ \vert \log  \vert 1-x \vert \vert}\Bigr\rceil -1;
3
\Bigr\},
\end{equation}
and  $r^\prime_{\Gamma}(x,\acc) = r_{\Gamma}(1+x,\acc)$.
For $x\in (1/2,1)$ there exists $\theta_1=\theta_1(x)\in (-1/2,1/2)$ such that
\begin{equation}
\label{Gamma-x>1/2}
\log \Gamma(x)
=
\gamma (1-x)  + \sum_{k=2}^{r_\Gamma(x,\acc)} \frac{\zeta(k)}{k}(1-x)^k
+\vert \theta_1 \vert 2^{-\acc}.
\end{equation}
For $x\in (0,1/2)$ we have that
there exists $\theta_2=\theta_2(x)\in (-1/2,1/2)$ such that
\begin{equation}  
\label{Gamma-x<1/2}
\log \Gamma(x)
=
- \log x   -\gamma x  + \sum_{k=2}^{r^\prime_\Gamma(x,\acc)} \frac{ \zeta(k)}{k} (-x) ^k
+\vert \theta_2 \vert 2^{-\acc}.
\end{equation}
\end{Proposition}

\begin{Proof}
From \eqref{coeff-gamma-def} we have that 
$0<c_{\Gamma}(k)<1$ for $k\ge 1$; recalling also \eqref{coeff-gamma-decreasing},
it is easy to see that in  this case \eqref{r-gen-estim} becomes \eqref{r-gamma-estim}.
For $x\in (1/2,1)$, we can compute $\log \Gamma(x)$
by specialising \eqref{f-x>1/2} thus  obtaining \eqref{Gamma-x>1/2}.
For $x\in (0,1/2)$, we recall $g_\Gamma(s,x) = \log x$
and $r^\prime_\Gamma(x,\acc) = r_\Gamma(1+x,\acc)$;
hence  \eqref{f-x<1/2}  becomes \eqref{Gamma-x<1/2}.
\end{Proof}

Moreover, in this case Proposition \ref{DIF-formulae-f} becomes
\begin{Proposition}[The reflection formulae for $\log \Gamma$]
Let $x\in (0,1)$, $x\ne 1/2$,
$\acc\in \N$, $\acc\ge 2$, $r_1(x,\acc) = r^\prime_\Gamma(x,\acc)/2 $
and $r_2(x,\acc) = r_\Gamma(x,\acc)/2 $.
There exists $\theta=\theta(x) \in (-1/2,1/2)$  such that
for $0<x  < 1/2$ we have
\begin{align*} 
\log\Gamma(x)&+\log\Gamma(1-x) 
= 
-  \log x + 
\sum_{\ell=1}^{r_1} \frac{\zeta(2\ell)}{\ell} x^{2\ell}
+ \vert \theta \vert 2^{-\acc}, 
\\
\log\Gamma(x)&-\log\Gamma(1-x) 
= 
-  \log x  - 2\gamma x 
-2 
\sum_{\ell=1}^{r_1} \frac{\zeta(2\ell+1)}{2\ell+1} x^{2\ell+1}
+ \vert \theta \vert 2^{-\acc}, 
\end{align*}
and  for $1/2<x<1$ we have
\begin{align*}
\log\Gamma(x)&+\log\Gamma(1-x) 
=
-  \log (1-x) + 
\sum_{\ell=1}^{r_2} \frac{\zeta(2\ell)}{\ell} (1-x)^{2\ell}
+ \vert \theta \vert 2^{-\acc},
\\ 
\log\Gamma(x)&-\log\Gamma(1-x) 
=
\log (1-x) +  2\gamma(1-x)
+2 \sum_{\ell=1}^{r_2} \frac{\zeta(2\ell+1)}{2\ell+1} (1-x)^{2\ell+1}
+ \vert \theta \vert 2^{-\acc}.
\end{align*}
\end{Proposition}
Similar formulae hold for the infinite series too.
We now see how to handle the tail of this case.

\subsection{The  Euler-Maclaurin formula for the tail of  $\log\Gamma$}
We now see how to evaluate the tail of the problem 
in \eqref{main-recursive-f},  
with $f(s,x) = \log\Gamma(x)$ and $g_{\Gamma}(s,x)= \log x$.
First of all, equation \eqref{horiz-shift-f} becomes 
\begin{align}
\log \Gamma(\lfloor x \rfloor + \{ x\}) 
-
\log \Gamma(\{ x\})   
\label{horiz-shift-Gamma}
&=
\log \Gamma( u + v )  
-
\log \Gamma(v) 
+
\sum_{j=0}^{t-1} \log (\{ x\}+ j  ) , 
\end{align}
where $t\in \N$,
$2\le t \le \lfloor x \rfloor$, $u:=\lfloor x \rfloor-t \in \N$ and $v := \{ x\}+t> t$. 
Now we can apply  the Euler-Maclaurin formula as in Lemma \ref{Euler-Maclaurin} to 
$\log \Gamma( u + v )   - \log \Gamma(v)$
in a situation, $v>t \ge 2$, in which we  have a small error term.  
Before applying Lemma \ref{Euler-Maclaurin} to $h(w)=\log w$, $w>0$, we need a 
couple of definitions. 
For  $v>0$, $u\in \N$, $u\ge 1$, $m\in\N$, $m\ge 1$, we define
\begin{equation}
\label{d-def}
d(u,v):= \Bigl(1+ \frac{u-1}{v}\Bigr)^{-1},
\end{equation}
\begin{equation}
\label{S-def}
S_{m}(u,v)
:=
\sum_{n=1}^m \frac{B_{2n}}{2n(2n-1)} 
\frac{d(u,v)^{2n-1} -1}{v^{2n-1}},
\end{equation}
where $B_{2n}$ are the even-index Bernoulli numbers
(we recall that, with the unique exception of $B_1=-1/2$,  the 
odd-index Bernoulli numbers are equal to $0$), and
\begin{equation}
\label{E-def}
E_{m}(u,v)
:=
\frac{1}{2m}
\int_{v}^{u + v-1}
B_{2m} ( \{t-v  \} )
\frac{\dx t} {  t^{2m}},
\end{equation}
where $B_{2m} ( u )$ are the even-index Bernoulli polynomials.
We remark that $S_{m}(v,1)   =0$ and $E_{m}(v,1) =0$.
We also define $E_{m}(0,u) :=0$ for every $u\in \N$, $u\ge 1$.
We will use the following 
\begin{Proposition}[The Euler-Maclaurin formula for $\log\Gamma$]
\label{LM-prop}
Let $v>0$, $u\in \N$, $u \ge 2$, $m\in\N$, $m\ge 1$. Then
\begin{equation}
\label{LM-main}
\log \Gamma  ( u + v ) - \log \Gamma  (v)
=
u\log v - u +1 -  \Bigl(u + v -\frac{1}{2} \Bigr)  
\log d(u,v)
+ S_{m}(u,v)  + E_{m}(u,v),
\end{equation}
where $\Gamma$ is  the Euler gamma-function, $d(u,v)$, $S_{m}(u,v)$, $E_{m}(u,v)$ are
respectively defined in \eqref{d-def}-\eqref{E-def}.
We further have that 
\begin{equation}
\label{LM-error-estim} 
\vert  E_{m}(u,v) \vert
<
2\Bigl(1+\frac{1}{4^{m+1}} \frac{2m+3}{2m+1} \Bigr)  
\frac{(2m)!}{(2\pi)^{2m+2}} 
\frac{1-d(u,v)^{2m+1}}{ v^{2m+1}}.
\end{equation}
Moreover,
for $v>0$, $u\in \{1,2\}$, we have 
%
\(
\log \Gamma  ( u + v ) - \log \Gamma  (v)
=
u\log v - (u-1)\log d(u,v).
\)
\end{Proposition}
\begin{Proof}
Let $v>0$, $u\in \N$, $u\ge 1$.
We use the Euler-Maclaurin summation formula on 
\(
\sum_{j=0}^{u-1} \log  (v + j  ). 
\)
Applying the first part of  Lemma \ref{Euler-Maclaurin} to $h(w)=\log w$, $a=v$
and $N=u-1$, we obtain
that
\begin{align} 
\notag
\sum_{j=0}^{u-1} \log (v + j )
&= 
\int_{v}^{u + v-1} \log w\ \dx w
+ 
\frac12  \bigl(\log (u + v-1 ) + \log v   \bigr)
\\&
\label{app-euler-maclaurin}
+
\sum_{n=1}^{m}  \frac{B_{2n}}{2n(2n-1)}
\Bigl( \frac{1}{(u + v-1)^{2n-1}}-  \frac{1}{v^{2n-1}} \Bigr)  
- \frac{1}{2m} 
\int_{v}^{u + v-1}  B_{2m} ( \{w-v\})\frac{\dx w}{w^{2m}},
\end{align}
since $(\log w)^{(\ell) }=  (-1)^{\ell-1}   w^{-\ell} (\ell-1)! $ 
for every $w>0$ and $\ell \ge 1$.
Recalling \eqref{d-def}, that $\int_{z_1}^{z_2} \log w \, \dx w= z_1 \log z_1 -z_1  -z_2 \log z_2 +z_2$ and
performing the needed computations in \eqref{app-euler-maclaurin},
the first part of Proposition \ref{LM-prop} follows.
Using the second part of Lemma \ref{Euler-Maclaurin},
we also have
\begin{align}
\notag
\Bigl\vert
\frac{1}{2m} 
\int_{v}^{u + v-1} B_{2m} ( \{w-v \}) \frac{ \dx w}{w^{2m}}
\Bigr\vert
&\le
\frac{\vert B_{2m+2} \vert (2m)!}{(2m+2)!} 
\Bigl \vert \frac{1}{(u + v-1)^{2m+1}}-  \frac{1}{v^{2m+1}}\Bigr \vert.
\end{align}
Recalling
\(
B_{2\ell}  = 2(-1)^{\ell-1} \frac{(2\ell)!}{(2\pi)^{2\ell}} \zeta (2\ell),
\)
where $\ell\in \N$, $\ell\ge 1$, 
by Lemma \ref{elementary-estim} one has
\(
\zeta (2\ell) 
<
1 + \frac{1}{4^{\ell}}  \frac{2\ell+1}{2\ell-1},
\)
so that 
\[
\vert B_{2m+2}  \vert < 2 \Bigl(1 + \frac{1}{4^{m+1}}
\frac{2m+3}{2m+1}  \Bigr) \frac{(2m+2)!}{(2\pi)^{2m+2}}
\]
and the second part of Proposition \ref{LM-prop} follows.
Letting $u\in \{1,2\}$, the third part of Proposition \ref{LM-prop} 
 follows from \eqref{d-def} and \eqref{difference-gamma}
since 
\(
\log \Gamma  (v+2) - \log \Gamma  (v)
=
2 \log v + \log(1+1/v).
\)
\end{Proof}
\subsubsection{Why we need the horizontal shift}
\label{why-hor-shift}
A well known characteristic of the Euler-Maclaurin summation formula,
see, e.g., Stoer-Bulirsch \cite[\S3.3-3.4]{StoerBulirsch2002}, 
is, for $u, v >0$ fixed, that 
$E_{m}(u,v)$ might not converge to $0$ as $m\to +\infty$.

In our application we  have $v=\{x\}<1$ and it is clear that for 
$v<1$, $E_{m}(\lfloor x\rfloor,\{x\})$, defined in \eqref{E-def}, diverges. 
But, after having used the horizontal shift as in \eqref{horiz-shift-Gamma},
 it is still possible to use Proposition \ref{LM-prop} to efficiently and
accurately evaluate $\log\Gamma( u + v )  - \log\Gamma(v)$, 
$u:=\lfloor x \rfloor-t \in \N$, $v := \{x\}+t \ge 2$, since
another well known characteristic of $E_{m}(u,v)$, $v>1$ and $u\in \N$ both fixed, 
is that  it  decreases for the first values of $m$ until it
reaches a minimal value for some $\widetilde{m}(v)$.
We can first determine  $\widetilde{m}(v)$
and evaluate $E_{\widetilde{m}}(u,v)$; if it is sufficiently small for our goals
we then proceed to evaluate the remaining part of \eqref{LM-main}.
Moreover, it is easy to see that the order of magnitude of the right hand side 
of \eqref{LM-error-estim}  depends weakly from large values of 
$u$; hence, for $v>1$ fixed,
we can find $\widetilde{m}$
such that 
\begin{equation}
\label{unifE-estim}
E_{m}(v)  
:=
\Bigl(1 + \frac{1}{4^{m+1}}  \frac{2m+3}{2m+1}\Bigr)  
\frac{(2m)!}{(2\pi)^{2m+2}}
\frac{2}{v^{2m+1}}
\end{equation}
is minimal  and then use $S_{\widetilde{m}}( u , v)$  and 
$E_{\widetilde{m}}(u,v)$ in \eqref{LM-main},
since $\vert E_{\widetilde{m}}(u,v) \vert \le   E_{\widetilde{m}}(v)$.
In practice, this procedure for $v>1$ works very well and its accuracy improves
as $v$ becomes larger.
\begin{Example}
\label{example-hor-shift}
As an example of this procedure, we consider the cases   $v=4,8,10,15$: 
a quick verification using a Pari/GP script gives that  
the optimal $\widetilde{m}$ for \eqref{unifE-estim}  are respectively $12,25,31,47$
and that  
$E_{12}(4)\approx 1.947867\cdot 10^{-12}$,
$E_{25}(8)\approx  1.664034\cdot 10^{-23}$,
$E_{31}(10)\approx5.192957\cdot 10^{-29}$,
$E_{47}(15)\approx 9.62651\cdot 10^{-43}$.
In Table \ref{EM-error-table} you can find the analysis for $v=4$.

\begin{table}[ht] 
\renewcommand{\arraystretch}{1.2}
\scalebox{1}{
\begin{tabular}{|c|c|}
\hline
$m$ & $E_{m} (4)$ \\ 
\hline
 $1$ & $<4.4278732056 \cdot10^{-5}$ \\\hline
 $2$ & $<7.7850220132 \cdot10^{-7}$ \\\hline
 $3$ & $<3.6364641750 \cdot10^{-8}$ \\\hline
 $4$ & $<3.2116629975 \cdot10^{-9}$ \\\hline
 $5$ & $<4.5719344226 \cdot10^{-10}$ \\\hline
 $6$ & $<9.5521140604 \cdot10^{-11}$ \\\hline
 $7$ & $<2.7521246133 \cdot10^{-11}$ \\\hline
 $8$ & $<1.0456683377 \cdot10^{-11}$ \\\hline
 $9$ & $<5.0656394025 \cdot10^{-12}$ \\\hline
 $10$ & $<3.0474585283 \cdot10^{-12}$ \\
\hline
\end{tabular} 
\begin{tabular}{|c|c|}
\hline
$m$ & $E_{m} (4)$ \\ 
\hline
 $11$ & $<2.2289481967 \cdot10^{-12}$ \\\hline
 $\bm{12}$ & $<1.9478670552 \cdot10^{-12}$ \\\hline
 $13$ & $<2.0044394624 \cdot10^{-12}$ \\\hline
 $14$ & $<2.3990263557 \cdot10^{-12}$ \\\hline
 $15$ & $<3.3042625797 \cdot10^{-12}$ \\\hline
 $16$ & $<5.1892728314 \cdot10^{-12}$ \\\hline
 $17$ & $<9.2176378732 \cdot10^{-12}$ \\\hline
 $18$ & $<1.8386982725 \cdot10^{-11}$ \\\hline
 $19$ & $<4.0927580309 \cdot10^{-11}$ \\\hline
 $20$ & $<1.0107900272 \cdot10^{-10}$ \\
\hline
\end{tabular} 
}
\caption{{\small Bounds for $E_{m}(4)$, $m=1,\dotsc,20$. The boldfaced $m$-value represents $\widetilde{m}$ for this case.}}
\vspace{-4mm}
\label{EM-error-table}
\end{table} 
This means that with the horizontal shifts $t=4,8,10,15$ we can evaluate
the tail of the $\log\Gamma$-function, namely $\log \Gamma  ( u + v ) -\log \Gamma  (v)$, 
using \eqref{LM-main} with an accuracy of, respectively, $32,64,80,128$ bits
uniformly for every $u \ge 2$.  
These are the values used to build Table \ref{Gamma-tail-accuracy-table}
below.
\end{Example}

Hence, by choosing $t$ large enough, and at the cost
of evaluating $t$ logarithms (or a logarithm of a product having $t$ factors;
a quantity that might be evaluated with a better computational cost if $t$ is not too large),  
we can obtain a much smaller estimate for
$E_{m}(u,v)$ than the one originally available for $E_{m}(\lfloor x \rfloor, \{x\})$.
For fixed precision computations we can use the parameters of 
Table \ref{Gamma-tail-accuracy-table}
we experimentally computed in Example \ref{example-hor-shift}. For multiprecision applications, 
we need in input a desired accuracy $\Delta\in (0,1)$ and, 
before launching the actual computation of $\log\Gamma$, we have to perform
the previously described evaluation of $\widetilde{m}(v)$ 
(and, if necessary, to determine a suitable large horizontal shift)
until we obtain that $E_{\widetilde{m}}(v)<\Delta$.
After this, we can evaluate the tail of $\log \Gamma$ using Proposition \ref{LM-prop}
being certain that the final result will have the required accuracy.
In the following we will refer to this procedure by saying that
\emph{the choice of the parameters must be performed at runtime}.

We finally remark that the maximal order of magnitude as  $x\to 0^+$ for $\log \Gamma(x)$
is $\log x$; so to get  its accurate evaluation 
one needs to work with at least 
$\lceil  \vert \log_2 \vert \log x \vert \vert \rceil$ binary digits,
where $\log_2(w)$ denotes the base-2 logarithm of $w$.
For $x \to +\infty$, Stirling's formula reveals that
we need at most $\lceil   \log_2 x +    \log_2 \log x  \rceil$ binary digits.

\begin{table}[ht] 
\renewcommand{\arraystretch}{1.2}
\scalebox{1}{
\begin{tabular}{|c|c|l|c|}
\hline
$t$ & $\widetilde{m}$ & \phantom{0123456}$E_{\widetilde{m}} (v)$ & binary precision\\ 
\hline
$4$ & $12$ & $\approx 1.947867\cdot 10^{-12}<2^{-33}$ & $32$\\ \hline
$8$ & $25$ & $\approx 1.664034\cdot 10^{-23}< 2^{-65}$ & $64$\\ \hline
$10$ & $31$ & $\approx 5.192957\cdot 10^{-29}< 2^{-81}$ & $80$\\ \hline
$15$ & $47$ & $\approx 9.626509\cdot 10^{-43} < 2^{-129}$ & $128$\\
\hline
\end{tabular} 
}
\caption{{\small Horizontal shift optimal error evaluation for $\log \Gamma$.}}
\vspace{-4mm}
\label{Gamma-tail-accuracy-table}
\end{table} 

\subsection{Computational costs and error terms for $\log \Gamma(x), x>0$}
Thanks to   Section \ref{comp-cost-f}, we can say that 
the cost of computing $\log \Gamma\{x\}$, $\{x\}\in (1/2,1)$
is  $\Odi{\acc}$ floating point products and $\Odi{\acc}$ floating point summations
with an accuracy of $\acc$ binary digits;
for $\{x\}\in (0,1/2)$ we have the same plus the cost of computing $\log \{x\}$.

We now evaluate the computational cost of the Euler-Maclaurin summation formula.
Let now  $\Delta\in(0,1)$ be fixed and $m :=m (\Delta,u,v)\in \N$.
We consider  the sum 
$S_{m}( u , v)$ as defined in   \eqref{S-def}. 
We can precompute and 
store the values of $B_{2\ell}/(2\ell(2\ell-1))$; for example
\begin{align*}
B_2 &= \frac16 , \quad
B_4 = -\frac{1}{30}, \quad
B_6 = \frac{1}{42}, \quad
B_8 =  -\frac{1}{30}, \quad
B_{10} = \frac{5}{66}, \quad
B_{12} = -\frac{691}{2730}, \quad
B_{14} = \frac76, \quad
\dotsc
\end{align*}
Then for a fixed $\ell$, to compute the summand in \eqref{S-def} 
we need four products and two sums since the powers of $v$ and $d(u,v)$
can be obtained exploiting a repeated product strategy.
So the cost of computing  $S_{m}( u , v)$ is
$4m$ products and $2m$  sums. 
To obtain a final result  having a $\Delta$-accuracy  using
\eqref{LM-error-estim}, we  need to first 
determine the smallest $t$ and the optimal $\widetilde{m}$ 
such that $E_{\widetilde{m}}(v) < \Delta$ using \eqref{unifE-estim}; then,
we choose the smallest $m \le \widetilde{m}$ 
such that $E_{m}(v)   < \Delta$.
To combine the previous analyses we have to choose $\Delta=2^{-\acc-1}$.
We also have to consider the cost of the horizontal shift procedure
which is just the cost of computing the sum of $t$ logarithms, or,
in alternative, the cost of a logarithm of a product having $t$ factors.

\section{The digamma function} 
\label{digamma-section}

First of all, we prove the following
\begin{Theorem}[The digamma-function is in $\scrF$]
$\psi(x)$, $x>0$, is in $\scrF$.
\end{Theorem}
\begin{Proof}
Recalling that  the well known Euler formula, see, e.g., Lagarias 
\cite[section 3]{Lagarias2013}, gives 
\begin{equation*}
\psi(z)  = - \gamma  - \sum_{k=1}^{\infty} \zeta(k+1) (1-z)^{k} ,
\end{equation*}
where $\zeta(\cdot)$ is the Riemann zeta-function, $\gamma$ 
is the Euler-Mascheroni constant and $z\in(0,2)$, by letting 
\begin{equation}
\label{coeff-psi-def}
c_{\psi}(k) = -\zeta(k+1)\ \textrm{for}\ k\ge 1, \quad  c_{\psi}(0) 
= -\gamma\quad \textrm{and}\quad 
C_{\psi}(k):= \vert c_{\psi}(k)\vert
\end{equation} 
we have that point \ref{series-f}) of  Definition \ref{set-def} 
holds. 
By Lemma \ref{elementary-estim} 
we have $1+ 2^{-w}<\zeta(w)<1+2^{1-w}$ for $w\ge 3$, so that
\begin{equation}
\label{coeff-psi-decreasing}
\Bigl\vert \frac{c_{\psi}(k+1)}{c_{\psi}(k)} \Bigr\vert  
=
\frac{C_{\psi}(k+1)}{C_{\psi}(k)} 
=
\frac{\zeta(k+2)}{\zeta(k+1)}
<
\frac{1+2^{-k-1}}{1+2^{-k-1}}
=
1
\quad 
\text{for}
\ k\ge 2.
\end{equation}
Hence point \ref{control-number-summands-f}) of Definition \ref{set-def} also holds
with  $\overline{k}_\psi =2$. 
Using the  functional equation  
\begin{equation}
\label{difference-psi}
\psi(1+z) = \psi(z)  +  \frac1z,
\quad z>0,
\end{equation}
we can say that point \ref{recursive-f}) of Definition \ref{set-def} holds 
with $g_{\psi}(s,z) = 1/z$, $z>0$.
This proves that $\psi(x)\in\scrF$.
\end{Proof}

The analysis on the number of summands that follows from Sections
\ref{comput-f-in-(0,1)}-\ref{comp-cost-f} gives  results similar to the ones already 
proved in Section 4 of \cite{Languasco2021}; the same happens for the 
reflection formulae we obtain by specialising Proposition \ref{DIF-formulae-f}.
Here we have just a slightly better estimate for $r_{\psi}(x,\acc)$.
Recalling that   $\psi(1/2)=-2\log 2 -\gamma$, we have the following
\begin{Proposition}[The number of summands for $\psi$]
Let $x\in (0,1)$, $x\ne 1/2$,
$\acc\in \N$, $\acc\ge 2$,
\begin{equation}
\label{r-psi-estim}
r_{\psi}(x,\acc) = 
\max
\Bigl\{
\Bigl\lceil \frac{(\acc+1) \log 2 
+ \vert \log (1-\vert 1-x\vert)\vert +0.2}{ \vert \log  \vert 1-x \vert \vert}\Bigr\rceil -1;
2
\Bigr\},
\end{equation}
and  $r^\prime_{\psi}(x,\acc) = r_{\psi}(1+x,\acc)$.
For $x\in (1/2,1)$ there exists $\theta_1=\theta_1(x)\in (-1/2,1/2)$ such that
\begin{equation}
\label{psi-x>1/2}
\psi(x)
= - \gamma  - \sum_{k=1}^{r_\psi(x,\acc)} \zeta(k+1) (1-x)^{k} 
 +\vert \theta_1 \vert 2^{-\acc}. 
\end{equation}
For $x\in (0,1/2)$ we have that
there exists $\theta_2=\theta_2(x)\in (-1/2,1/2)$ such that
\begin{equation}  
\label{psi-x<1/2}
\psi(x)
=
- \frac{1}{x}  -\gamma   - \sum_{k=1}^{r^\prime_\psi(x,n)}  \zeta(k+1) (-x)^{k}
+\vert \theta_2 \vert 2^{-\acc}.
\end{equation}
\end{Proposition}
\begin{Proof}
Recalling  \eqref{coeff-psi-decreasing} and remarking that  $\log(\zeta(3))<0.185 < 0.2$,
it is easy to see that in  this case \eqref{r-gen-estim} becomes \eqref{r-psi-estim}.
For $x\in (1/2,1)$, we can compute $\psi(x)$
by specialising \eqref{f-x>1/2} thus  obtaining \eqref{psi-x>1/2}.
For $x\in (0,1/2)$, we recall $g_\psi(s,x) = 1/ x$
and $r^\prime_\psi(x,\acc) = r_\psi(1+x,\acc)$;
hence  \eqref{f-x<1/2}  becomes \eqref{psi-x<1/2}.
\end{Proof}

Moreover, recalling \eqref{coeff-psi-def}, in this case Proposition \ref{DIF-formulae-f} becomes
\begin{Proposition}[The reflection formulae for $\psi$]
Let $x\in (0,1)$, $x\ne 1/2$,
$\acc\in \N$, $\acc\ge 2$, $r_1(x,\acc) = r^\prime_\psi(x,\acc)/2 $
and $r_2(x,\acc) = r_\psi(x,\acc)/2 $.
There exists $\theta=\theta(x) \in (-1/2,1/2)$  such that
for $0<x  < 1/2$ we have
\begin{align*}  
\psi(x)&+\psi(1-x) 
=
-2\gamma
-\frac{1}{x}  - 
2\sum_{\ell=1}^{r_1} \zeta(2\ell+1) x^{2\ell} 
+ \vert \theta \vert 2^{-\acc},\\
\psi(x)&-\psi(1-x) 
=
- \frac{1}{x}  +
2\sum_{\ell=1}^{r_1} \zeta(2\ell) x^{2\ell-1} 
+ \vert \theta \vert 2^{-\acc} ,
\end{align*}
and  for $1/2<x<1$ we have
\begin{align*}
\psi(x)&+\psi(1-x) 
=
-2\gamma
-\frac{1}{1-x}  
-
2\sum_{\ell=1}^{r_2}  \zeta(2\ell+1) (1-x)^{2\ell} 
+ \vert \theta \vert 2^{-\acc},
\\
\psi(x)&-\psi(1-x) 
=
\frac{1}{1-x}  -
2\sum_{\ell=1}^{r_2}  \zeta(2\ell) (1-x)^{2\ell-1} 
+ \vert \theta \vert 2^{-\acc} .
\end{align*}
\end{Proposition}
Similar formulae hold for the infinite series too.
We now see how to handle the tail of this case.

\subsection{The  Euler-Maclaurin formula for the tail of   $\psi(x)$, $x>0$}
\label{psi-tail}
We now see how to evaluate the tail of the problem in \eqref{main-recursive-f},  
with $f(s,x) =\psi(x)$ and $g_{\psi}(s,x)= 1/ x$, $x>0$.
First of all, equation \eqref{horiz-shift-f} becomes 
\begin{align*}
\psi(\lfloor x \rfloor + \{ x\}) 
-
\psi(\{ x\})  
&= 
\psi(u + v )  
-
\psi(v) 
+
\sum_{j=0}^{t-1} \frac{1}{\{ x\}+ j},
\end{align*}
where $t\in \N$,
$2\le t \le \lfloor x \rfloor$, $u:=\lfloor x \rfloor-t \in \N$ and $v := \{ x\}+t> t$. 
Now we can apply  the Euler-Maclaurin formula to $\psi(u + v )   - \psi(v)  $
in a situation, $v>t \ge 2$, in which we  have a small error term.  

The computation of $\sum_{j=0}^{u-1} (v+j)^{-1}$
can be performed using Lemma \ref{Euler-Maclaurin}  with $g_f(s,z) = z^{-1}$.
In this way we get that
\begin{Proposition}[The Euler-Maclaurin formula for $\psi$]
\label{LM-digamma-EM}
Let $v>0$, $u\in \N$, $u \ge 2$, $m\in\N$, $m\ge 1$. Then
\begin{align}
\notag
\psi( u + v ) - \psi(v)
&
= 
- \log d(u,v) 
+  
\frac{1}{2v}\bigl(1 + d(u,v)\bigr)
+
\sum_{n=1}^{m} \frac{B_{2n}}{2n} 
\frac{1 - d(u,v)^{2n}}{v^{2n}} 
\\&\hskip1cm
\label{digamma-EM}
-
\int_{v}^{u + v-1} 
B_{2m} (\{w-v \})  \frac{\dx w}{w^{2m+1}},
\end{align}
where 
$B_{2n}$ are the even-index \emph{Bernoulli numbers},
$B_{2n}(u)$ are the even-index \emph{Bernoulli polynomials} and $d(x,y)$ is
defined in \eqref{d-def}.
We further have that 
\begin{equation}
\label{digamma-EM-error}
\Bigl \vert 
\int_{v}^{u + v-1} 
B_{2m} ( \{w-v \})  \frac{\dx w}{w^{2m+1}}
\Bigr \vert
<
2\Bigl(1 + \frac{1}{4^{m+1}}  \frac{2m+3}{2m+1}\Bigr)  
\frac{(2m+1)!}{(2\pi)^{2m+2}} 
\frac{1 - d(u,v)^{2m+2}}{v^{2m+2}}.
\end{equation}
\end{Proposition}
The proof of  Proposition \ref{LM-digamma-EM} is completely analogous to the one of
Proposition \ref{LM-prop}. 

\subsection{Computational costs and error terms for $\psi(x), x>0$}
Thanks to  Section \ref{comp-cost-f}, we can say that 
the cost of computing $\psi\{x\}$, $\{x\}\in (1/2,1)$
is  $\Odi{\acc}$ floating point products and $\Odi{\acc}$ floating point sums
with an accuracy of $\acc$ binary digits;
for $\{x\}\in (0,1/2)$ we have the same plus the cost of computing $1/\{x\}$.

Let now  $\Delta\in(0,1)$ be fixed and $m :=m (\Delta,u,v)\in \N$.
We consider  the sum   in   \eqref{digamma-EM}. 
In this case too the key point is that   we can
store the values of $B_{2\ell}/(2\ell)$.
Then
for a fixed $n$, to compute the summand in \eqref{S-def} 
we need four products and two sums since the powers of $v$ and $d(u,v)$
can be obtained exploiting a repeated product strategy.
So the cost of computing  the sum   in   \eqref{digamma-EM} is
$4m$ products and $2m$  sums.  
To obtain
a final result  having a $\Delta$-accuracy  using
\eqref{LM-error-estim}, we  need to first 
determine the smallest $v$ and the optimal $\widetilde{m}(v)$ such that
\begin{equation}
\label{error-EM-psi}
E^{\psi}_{m}(v) := 
\Bigl(1 + \frac{1}{4^{m+1}}  \frac{2m+3}{2m+1}\Bigr)  
\frac{(2m+1)!}{(2\pi)^{2m+2}}
\frac{2}{v^{2m+2}}
\end{equation}
is minimal and less than $\Delta$ (remark that the right hand side of \eqref{digamma-EM-error} 
is  $ \le E^{\psi}_{m}(v)$). Then
we choose the smallest $m(v, \Delta) \le \widetilde{m}(v)$ 
such that $E^{\psi}_{m}(v)   < \Delta$.
In this procedure it might be necessary to use a suitable horizontal
shift $t$;
it is enough to  replicate the discussion in Section \ref{why-hor-shift}
replacing $E_{m}(v)$ with $E^{\psi}_{m}(v)$.
To combine the previous analyses we have to choose $\Delta=2^{-\acc-1}$.
We also have to consider the cost of the horizontal shift procedure
which is just the cost of computing the sum of $t$ fractions.
For fixed precision computations we can use the parameters of Table 
\ref{psi-tail-accuracy-table} we experimentally computed in a totally analogous way as we did in 
Example \ref{example-hor-shift}; 
for multiprecision applications, such a choice
has to be performed at runtime, see the discussion in Section \ref{why-hor-shift}.

We finally remark that the maximal order of magnitude  as  $x\to 0^+$ for $\psi(x)$
is $1/x$; so to get  its accurate evaluation 
one needs to work with at least  $\lceil  \vert \log_2 x \vert \rceil$ binary digits.
For $x \to +\infty$, the well known asymptotic properties of $\psi(x)$ reveal that
we need at most $\lceil \log_2 \log x  \rceil$ binary digits.
\begin{table}[ht] 
\renewcommand{\arraystretch}{1.3}
\scalebox{1}{
\begin{tabular}{|c|c|l|c|}
\hline
$t$ & $\widetilde{m}$ & \phantom{0123456}$E^{\psi}_{\widetilde{m}} (v)$ & binary precision\\ \hline
$4$ & $12$ & $\approx 1.21741\cdot 10^{-11}<2^{-33}$ & $32$\\ \hline
$8$ & $24$ & $\approx 1.05109\cdot 10^{-22}< 2^{-65}$ & $64$\\ \hline
$10$ & $31$ & $\approx 3.27156\cdot 10^{-28}< 2^{-81}$ & $80$\\ \hline
$15$ & $46$ & $\approx 6.06445\cdot 10^{-42} < 2^{-129}$ & $128$\\
\hline
\end{tabular} 
}
\caption{{\small Horizontal shift optimal error evaluation for $\psi$.}}
\vspace{-4mm}
\label{psi-tail-accuracy-table}
\end{table}

\section {The Hurwitz zeta-function}
\label{Hurwitz-zeta-section}

As we already did for $\log\Gamma$ and $\psi$, the first result we prove is the following
\begin{Theorem}[The Hurwitz zeta-function is in $\scrF$]
The Hurwitz zeta-function
$\zeta(s,x)$, $x>0$, $s>1$ being fixed, is in $\scrF$.
\end{Theorem}
\begin{Proof}
Letting $s> 1$, $z\in(0,2)$ and 
using $\frac{\partial}{\partial z} \zeta(s, z) = -z \zeta(s+1, z)$,
it is not hard to prove that
\begin{equation}
\label{Euler-Hurwitz} 
\zeta(s, z) =      
\sum_{k=0}^{\infty}
\frac{\Gamma(k+s)}{(k!) \Gamma(s)}  \zeta(k+s) (1-z)^{k}. 
\end{equation}
For $s>1$ fixed, unfortunately, $\Gamma(k+s)\zeta(k+s)/[(k!) \Gamma(s)] $ 
is not a decreasing sequence\footnote{Formula \eqref{Euler-Hurwitz}, for $s\in \N$, 
$s\ge 2$ (see also Section \ref{polygamma-sec}),
is used in the \texttt{boost} software library, see
\url{https://www.boost.org/doc/libs/1_81_0/libs/math/doc/html/math_toolkit/sf_gamma/polygamma.html}.
Since the coefficients in  \eqref{Euler-Hurwitz} are not a decreasing sequence, it is not clear how they
were able to control the needed number of summands required to compute the values of these functions 
with a fixed accuracy $\Delta$.} 
in $k$. To overcome this problem, we isolate the power series centered at $1$ of $z^{-s} -1$ in \eqref{Euler-Hurwitz} thus
obtaining
\begin{equation}
\label{enlarged-radius-Hurwitz}
\zeta(s,z)
=
\zeta(s) -1  + \frac{1}{z^s} 
+
\sum_{k=1}^{\infty}
\frac{\Gamma(k+s)}{(k!) \Gamma(s)} (\zeta(k+s) -1)(1-z)^{k}. 
\end{equation}
Using Lemma \ref{elementary-estim},
the series in \eqref{enlarged-radius-Hurwitz} absolutely converges for 
every $z\in(-1,3)$. This proves 
that $\zeta(s,z)$ satisfies  point \ref{series-f}) of Definition \ref{set-def}
with $c_{\zeta_{H}}(s,0) = \zeta(s) -1  + z^{-s} $, 
\begin{equation}
\label{coeff-Hurwitz-def}
c_{\zeta_{H}}(s,k) 
=  
\frac{\Gamma(k+s)}{(k!) \Gamma(s)} (\zeta(k+s)-1)
=
\frac{\zeta(k+s)-1}{k B(s,k)}
\quad \textrm{and} \quad
C_{\zeta_{H}}(s,k)  :=c_{\zeta_{H}}(s,k)
\end{equation}
for $k\ge 1$ and  $s > 1$, where \(
B(s,k) := \Gamma(s)\Gamma(k)/ \Gamma(s+k)
\)
is the Euler beta-function. 
By Lemma \ref{elementary-estim}  and  the functional equation
$\Gamma(w+1) = w \Gamma(w)$,  $w>0$,  we obtain
\begin{equation}
\label{coeff-hurwitz-decreasing}
\frac{c_{\zeta_{H}}(s,k+1)}{c_{\zeta_{H}}(s,k)} 
=
\frac{\zeta(k+s+ 1)-1}{\zeta(k+s)-1}
\frac{k B(s,k)}{(k+1) B(s,k+1)}
<
\frac{1}{2}
\frac{k+s+2}{k+s}
\frac{k+s}{k+1}
=
\frac{1}{2}
\frac{k+s+2}{k+1}
\le
1
\end{equation}
for $k\ge s$.
Hence point \ref{control-number-summands-f}) of Definition \ref{set-def} holds 
with $k_{\zeta_{H}}(s) = \lceil s\rceil$. 
In this case the functional equation is
\begin{equation}
\label{difference-Hurwitz}
\zeta(s,z+1) = \zeta(s,z) - z^{-s},
\end{equation}
for every fixed $s>1$ and $z>0$;
hence  point \ref{recursive-f}) of Definition \ref{set-def} holds with 
$g_{\zeta_H}(s,z) =  -z^{-s}$.
Summarising, $\zeta(s,x)$, $s>1$ being fixed, $x>0$, is in $\scrF$.
\end{Proof}

We now perform the analysis on the number of summands that follows from Sections
\ref{comput-f-in-(0,1)}-\ref{comp-cost-f}. 
Remarking that $\zeta (s,1/2) = (2^{s}-1) \zeta(s)$, we have the following
\begin{Proposition}[The number of summands for $\zeta(s,x)$]
Let $s>1$ be fixed. Let further  $x\in (0,1)$, $x\ne 1/2$,
$\acc\in \N$, $\acc\ge 2$,
\begin{equation}
\label{r-Hurwitz-estim}
r_{\zeta_{H}}(s,x,\acc) = 
\max
\Bigl\{
\Bigl\lceil \frac{(\acc+1) \log 2 + \vert \log ( 1-\vert 1-x\vert ) \vert 
+ \vert \log c_{\zeta_{H}}(s,\lceil s\rceil+1)\vert}
{ \vert \log  \vert 1-x \vert \vert}\Bigr\rceil -1;
\lceil s\rceil 
\Bigr\},
\end{equation}
and  $r^\prime_{\zeta_{H}}(s,x,\acc) = r_{\zeta_{H}}(s,1+x,\acc)$.
For $x\in (1/2,1)$ there exists $\theta_1=\theta_1(s,x)\in (-1/2,1/2)$ such that
\begin{equation} 
\label{zeta-Hurwitz-x>1/2}
\zeta (s,x)  
=  
\zeta(s) -1 + \frac{1}{x^s} 
+
\sum_{k=1}^{r_{\zeta_{H}}(s,x,\acc)}  
\frac{\zeta(k+s)-1}{k B(s,k)}  (1-x)^{k}
+\vert \theta_1 \vert 2^{-\acc}.
\end{equation}
For $x\in (0,1/2)$ we have that
there exists $\theta_2=\theta_2(s,x)\in (-1/2,1/2)$ such that
\begin{equation} 
\label{zeta-Hurwitz-x<1/2}
\zeta (s,x)  
= 
\zeta(s) -1 + \frac{1}{(1+x)^s} +
\frac{1}{x^{s}} 
+  \sum_{k=1}^{r^\prime_{\zeta_{H}}(s,x,\acc)} 
\frac{\zeta(k+s)-1}{k B(s,k)}  (-x)^{k} 
+\vert \theta_2 \vert 2^{-\acc}.
\end{equation}
\end{Proposition}
\begin{Proof}
Recalling  \eqref{coeff-hurwitz-decreasing},
it is easy to see that in  this case \eqref{r-gen-estim} becomes \eqref{r-Hurwitz-estim}.
For $x\in (1/2,1)$, we can compute $\zeta(s,x)$
by specialising \eqref{f-x>1/2} thus  obtaining \eqref{zeta-Hurwitz-x>1/2}.
For $x\in (0,1/2)$, thanks to \eqref{difference-Hurwitz} we have $g_{\zeta_H}(s,x) =  -1/ x^s$
and $r^\prime_{\zeta_H}(s,x,\acc) = r_{\zeta_H}(s,1+x,\acc)$;
hence  \eqref{f-x<1/2}  becomes \eqref{zeta-Hurwitz-x<1/2}.
\end{Proof} 
The reflection formulae are collected in the following 
\begin{Proposition}[The reflection formulae for $\zeta(s,x)$]
\label{DIF-formulae-Hurwitz}
Let $s>1$ be fixed.
Let futher $x\in (0,1)$, $x\ne 1/2$,
$\acc\in \N$, $\acc\ge 2$, $r_1(s,x,\acc)=r^\prime_{\zeta_{H}}(s,x,\acc)/2$
and $r_2(s,x,\acc)$ $=r_{\zeta_{H}}(s,x,\acc)/2$.
Recalling \eqref{coeff-Hurwitz-def}, there exists 
$\theta=\theta(s,x)\in (-1/2,1/2)$  such that 
for $0<x <1/2$ we obtain
\begin{align}
\notag
\zeta(s,x) +\zeta(s,1-x)
&=
2 \zeta(s) -2 + 
\frac{1}{x^s}
+ \frac{1}{(1-x)^s}
+ \frac{1}{(1+x)^s}
+ 
2\sum_{\ell=1}^{r_1}  \frac{\zeta(2\ell+s)-1}{2\ell B(s,2\ell)} x^{2\ell}  
+ \vert \theta \vert 2^{-\acc},
\\
\zeta(s,x) -\zeta(s,1-x)
&=
\frac{1}{x^s} 
-  \frac{1}{(1-x)^s}
+  \frac{1}{(1+x)^s} 
\label{zeta-Hurwitz-DIF-odd-x<1/2}
-
2\sum_{\ell=1}^{r_1}
\frac{\zeta(2\ell-1+s)-1}{(2\ell-1) B(s,2\ell-1)}  x^{2\ell-1}
+ \vert \theta \vert 2^{-\acc} ,
\end{align}
and for $1/2<x <1$ we have
\begin{align*}
\zeta(s,x) +\zeta(s,1-x)
\notag
&=
2 \zeta(s) -2 
+ \frac{1}{x^s}  + \frac{1}{(1-x)^s}  + \frac{1}{(2-x)^s}
+
2\sum_{\ell=1}^{r_2}  \frac{\zeta(2\ell+s)-1}{2\ell B(s,2\ell)}(1-x)^{2\ell}  
+ \vert \theta \vert 2^{-\acc} ,
\\
\notag
\zeta(s,x) -\zeta(s,1-x) 
&=
\frac{1}{x^s}   -\frac{1}{(1-x)^s} -\frac{1}{(2-x)^s} 
+
2\sum_{\ell=1}^{r_2} \frac{\zeta(2\ell-1+s)-1}{(2\ell-1)B(s,2\ell-1)}(1-x)^{2\ell-1}  
+ \vert \theta \vert 2^{-\acc} .
\end{align*}
\end{Proposition}
Similar formulae hold for the infinite series too.

\begin{Remark} [Coefficients computations] 
We  remark that 
$b_{s,k} := \Gamma(k+s)/[(k!) \Gamma(s)]$
can be easily computed using a repeated product strategy since 
$b_{s,0} = 1$, $b_{s,1} = s$, $b_{s,2} = \frac{s+1}{2}b_{s,1}$,
\begin{equation}
\label{recursive-coeff-hurwitz}
b_{s,k+1}
= 
\frac{\Gamma(k+1+s)}{(k+1)! \ \Gamma(s)}
=
\frac{s+k}{k+1}  b_{s,k}
=
\Bigl(1+  \frac{s-1}{k+1}\Bigr)  b_{s,k},
\end{equation}
in which we used that $\Gamma(w+1) = w \Gamma(w)$, $w>0$.
\end{Remark}

\subsection{The Euler-Maclaurin formula for the tail of $\zeta(s,x)$, $s>1$ being fixed, $x>0$}
\label{Hurwitz-tail}
We now see how to evaluate the tail of the problem in \eqref{main-recursive-f},  
with $f(s,x) =\zeta(s,x)$ and $g_{\zeta_H}(s,x)=  -x^{-s}$, $s>1$ being fixed and $x>0$.
First of all, equation \eqref{horiz-shift-f} becomes 
\begin{align}
\notag
\zeta(s,\lfloor x \rfloor + \{ x\}) 
-
\zeta(s,\{ x\})   
\notag
&=
\zeta(s, u + v )  
-
\zeta(s,v) 
-
\sum_{j=0}^{t-1} (\{ x\}+j)^{-s} ,
\end{align}
where $t=t(s)\in \N$,
$2\le t \le \lfloor x \rfloor$, $u:=\lfloor x \rfloor-t \in \N$ and $v := \{ x\}+t> t$. 
Now we can apply  the Euler-Maclaurin formula to 
$\zeta(s, u + v )   - \zeta(s,v)  $
in a situation, $v>t \ge 2$, in which, as we will see in the next subsection, 
we have a small error term.  
The computation of $ \sum_{j=0}^{u-1} (v+j)^{-s}$
can be performed using Lemma \ref{Euler-Maclaurin}  with $g_{\zeta_H}(s,z) = -z^{-s}$.
In this way we get the following
\begin{Proposition}[The Euler-Maclaurin formula for $\zeta(s,x)$]
\label{LM-Hurwitz-EM}
Let $s>1$ be fixed. Let further $v>0$, $u\in \N$, $u \ge 2$, $m\in\N$, $m\ge 1$. Then
\begin{align}
\notag
\zeta(s, u + v) - \zeta(s,v)
&=
 \frac{d(u,v)^{s-1} - 1}{(s-1)\ v^{s-1}} 
-
\frac{1+d(u,v)^{s}}{2v^s}  
-
\sum_{n=1}^{m} \frac{B_{2n}}{2n}   \Bigl(\prod_{j=0}^{2n-2} \frac{s+j}{j+1} \Bigr) 
\frac{1-d(u,v)^{2n+s-1} }{v^{2n+s-1}}
\\&
\label{Hurwitz-EM}
\hskip1.5cm
-
\frac{\prod_{j=0}^{2m-1} (s+j)}{(2m)!}
\int_{v}^{u + v-1} B_{2m} ( \{w- v \} )  \frac{\dx w}{w^{2m+s}},
\end{align}
where 
$B_{2n}$ are the even-index \emph{Bernoulli numbers},
$B_{2n}(u)$ are the even-index \emph{Bernoulli polynomials} and $d(x,y)$ is
defined in \eqref{d-def}.
We further have that 
\begin{equation}
\label{Hurwitz-EM-error}
\Bigl \vert  
\frac{1}{(2m)!}
\int_{v}^{u + v-1} B_{2m} ( \{w- v \} )  \frac{\dx w}{w^{2m+s}}
\Bigr \vert
<
\Bigl(1 + \frac{1}{4^{m+1}} \frac{2m+3}{2m+1}\Bigr)  
\frac{1-d(u,v)^{2m+s+1}}{\pi(2\pi)^{2m+1}\ v^{2m+s+1}}. 
\end{equation}
\end{Proposition}
The proof of  Proposition \ref{LM-Hurwitz-EM} is completely analogous to the one of
Proposition \ref{LM-prop}. 

\subsection{Computational costs and error terms for $\zeta(s,x)$, $s > 1$ being fixed, $x>0$}
%
Thanks to  Section \ref{comp-cost-f}, we can say that 
the cost of computing
the sum in \eqref{zeta-Hurwitz-x>1/2} 
is $\Odi{\acc}$ floating point products and $\Odi{\acc}$ floating point sums
with an accuracy of $\acc$ binary digits.
To obtain the values  $\zeta(s,\{x\})$, $\{x\}\in (1/2,1)$, $s>1$ being fixed,
we also need to further add the cost of computing $\{x\}^{-s}$.
For $\{x\}\in (0,1/2)$ we have the same plus the cost of computing $(1+\{x\})^{-s}$.

For estimating the computational cost of using the Euler-Maclaurin formula,
we can repeat the same argument already used for $\log\Gamma$ and $\psi$.
Let now  $\Delta\in(0,1)$ be fixed and $m :=m(s,\Delta,u,v)\in \N$.
We consider  the sum   in   \eqref{Hurwitz-EM}. 
In this case too the key point is that   we can
store the values of $B_{2\ell}/(2\ell)$.
Hence, for a fixed $n$, to compute the summand in \eqref{Hurwitz-EM} 
we need five products and three sums, since  $b_{s,n}=\prod_{j=0}^{n-1} (s+j)/(j+1)$ can be obtained 
as in \eqref{recursive-coeff-hurwitz} and  the powers of $v$ and $d(u,v)$
can be obtained exploiting a repeated product strategy.
So the cost of computing  the sum   in   \eqref{Hurwitz-EM} is
$5m$ products and $3m$  sums.  
To obtain a final result  having a $\Delta$-accuracy  using
\eqref{Hurwitz-EM}, we  need to first 
remark that, using  \eqref{Hurwitz-EM-error}, 
one has that 
\[\Bigl \vert 
\frac{\prod_{j=0}^{2m-1} (s+j)}{(2m)!}
\int_{v}^{u + v-1} B_{2m} ( \{w- v \} )  \frac{\dx w}{w^{2m+s}}
\Bigr \vert
\le
\Bigl(1 + \frac{1}{4^{m+1}}  \frac{2m+3}{2m+1}\Bigr) 
\frac{\prod_{j=0}^{2m} (s+j)}{\pi (2\pi)^{2m+1}\ v^{2m+s+1}}. 
\]
So it will be enough
determine the smallest $v$ and the optimal $\widetilde{m}(s,v)$ such that
\begin{equation}
\label{error-EM-Hurwitz}
E^{\zeta_{H}}_{m}(s,v) := 
\Bigl(1 + \frac{1}{4^{m+1}}  \frac{2m+3}{2m+1}\Bigr) 
\frac{\prod_{j=0}^{2m} (s+j)}{\pi (2\pi)^{2m+1}\ v^{2m+s+1}} 
\end{equation}
is minimal and less than $\Delta$;  then,
we choose the smallest $m(s,v,\Delta) \le \widetilde{m}(s,v)$ 
such that $E^{\zeta_{H}}_{m}(s,v)  < \Delta$.
In this procedure it might be necessary to use a suitable horizontal
shift $t=t(s)$.
To combine the previous analyses we have to choose $\Delta=2^{-\acc-1}$.
We also have to consider the cost of the horizontal shift procedure
(which is  the cost of computing the sum of $t$ powers with real exponents).
The dependence from $s$ of $E^{\zeta_{H}}_{m} (s,v)$   prevents
us to \emph{a priori} establish what is the optimal level for the horizontal shift parameter as we did 
for $\log \Gamma$ and $\psi$;
hence this has to be computed at runtime  for both the fixed precision and the multiprecision
cases (see also the discussion in Section \ref{why-hor-shift}).

We finally remark that  the maximal order of magnitude as $x\to 0^+$, $s>1$ being fixed, for $\zeta(s,x)$
is $x^{-s}$; so to get  its accurate evaluation 
one needs to work with at least   $\lceil s \vert \log_2 x \vert \rceil$ binary digits.
For $x \to +\infty$, $s>1$ being fixed, the well known asymptotic properties of $\zeta(s,x)$ reveal that
we need at least $\lceil (s-1) \log_2  x  \rceil$ binary digits.

\subsection{The polygamma functions}
\label{polygamma-sec}
Let $w \in \N$,  $w\ge 1$.
The polygamma function $\psi^{(w)}(z)$ is defined as the $w$-th
derivative of the digamma function $\psi(z)$, $z>0$. 
Instead of writing an \emph{ad hoc} treatment for such functions, we use the fact that
\[
\psi^{(w)}(z) =
(-1)^{w-1} (w!) \zeta(w+1, z),
\]
where $\zeta(s,z)$ is the Hurwitz zeta-function.
This  implies that $\psi^{(w)}(z)\in \scrF$.
We also recall that  
$\psi^{(w)}(1/2) = (-1)^{w-1} (w!) (2^{w+1}-1) \zeta(w+1)$
and 
$\psi^{(w)}(1) = (-1)^{w-1} (w!) \zeta(w+1)$.
We clearly have that $c_{\psi^{(w)}}(k) := (-1)^{w-1} (w!)c_{\zeta_{H}}(k)$,
$C_{\psi^{(w)}}(k) := (w!)c_{\zeta_{H}}(k)$,
and 
\[
\frac{C_{\psi^{(w)}}(k+1)}{C_{\psi^{(w)}}(k)}  
=
\frac{c_{\zeta_{H}}(w+1,k+1)}{c_{\zeta_{H}}(w+1,k)}.
\]

Hence we can argue as we did in Section \ref{Hurwitz-zeta-section} 
for $\zeta(s, z)$, $z\in (0,1)$, since it is enough to change any occurrence of $c_{\zeta_{H}}(s,k)$ with 
$(-1)^{w-1}(w!)c_{\zeta_{H}}(w+1,k)$.
But, since $C_{\psi^{(w)}}(k) = (w!)c_{\zeta_{H}}(k)$, the value $\log(w!)$ has to be inserted at the 
numerator of \eqref{r-Hurwitz-estim} thus enlarging the number of summands needed
to obtain an accuracy of $\acc$ binary digits.

Moreover, for $z>1$, in adapting the Euler-Maclaurin formula 
application of Section \ref{Hurwitz-tail},
we need to work with 
\[
E^{\psi^{(w)}}_{m}\!\!(v) 
:= 
(w!) E^{\zeta_{H}}_{m}(w+1,v) 
\]
where $E^{\zeta_{H}}_{m}(s,v) $ is defined in \eqref{error-EM-Hurwitz}.
This means that we  also need to consider the $w!$-term present in $E^{\psi^{(w)}}_{m}\!\!(v)$
to be able to perform the computation with the desired accuracy $\Delta\in(0,1)$.
As a consequence, we will either have to increase the value of the horizontal shift $t(w+1)$
or, if possible, to choose a value of $m(w+1,\Delta,v)$ larger than the one needed to compute
$\zeta(w+1, z)$ with a $\Delta$-accuracy.

\section {The Hurwitz zeta-function: first partial derivative}
\label{Hurwitz-prime-section}

Now we show that $\zeta^{\prime}(s,x) = \frac{\partial\zeta}{\partial s}(s,x)$, $x>0$, $s>1$,  is in $\scrF$. 
\begin{Theorem}[The first partial derivative of the Hurwitz zeta-function is in $\scrF$]
The first partial derivative of the Hurwitz zeta-function,
$\zeta^{\prime}(s,x)= \frac{\partial\zeta}{\partial s}(s,x)$, $x>0$, $s>1$ being fixed, is in $\scrF$.
\end{Theorem}
\begin{Proof}
For $\zeta^{\prime}(s,z)$, $s>1$ being fixed and $z>0$, we proceed as follows.
Differentiating \eqref{enlarged-radius-Hurwitz} we obtain, 
for $s> 1$, $z\in(0,2)$, that
\begin{equation}
\label{Euler-Hurwitz-deriv}  
\zeta^{\prime}(s, z)  
=
\zeta^{\prime}(s)
-
\frac{\log z}{z^s}
+
\sum_{k=1}^{\infty} 
c_{\zeta^{\prime}_{H}}(s,k) 
(1-z)^{k},
\end{equation}
where 
\begin{equation}
\label{coeff-hurwitz-prime-def}
c_{\zeta^{\prime}_{H}}(s,k) : = \frac{\alpha(s,k)} {B(s,k)}, 
\quad 
\alpha(s,k): =\frac{(\zeta(k+s)-1) (\psi(k+s)- \psi(s)) + \zeta^{\prime}(k+s)}{k}
\end{equation}
and $B(s,k)$ is the Euler beta-function.
Thanks to the well known results about the radius of convergence of power series,
the series in \eqref{Euler-Hurwitz-deriv} trivially converges for every $z\in(-1,3)$. 
This proves that point \ref{series-f}) of Definition \ref{set-def} holds.
By Lemma \ref{elementary-estim}, for  $k\ge s+1$, $s>1$, and $k\ge 3$ we have 
\begin{align}
\notag
\vert \alpha(s,k) \vert 
&<
\frac{1}{k2^{k+s}}\Bigl(\log(k+s) \frac{k+s+1}{k+s-1}+  1.43\Bigr)
\\&
\label{alpha-estim}
\le
\frac{1}{2^{k+s}} 
\Bigl(\max\Bigl(\frac{4\cdot \log 5}{9};\frac{\log(2s+1)}{s}\Bigr)+\frac{1.43}{3}\Bigr)
< \frac{1.64}{2^{k+s}}.
\end{align}
Hence for $k\ge \max(s+1;3)$ and $s>1$ we have 
$\vert c_{\zeta^{\prime}_{H}}(s,k)\vert <1.64/(2^{k+s}B(s,k))=: C_{\zeta^{\prime}_{H}}(s,k)$.
Since
\[
\frac{C_{\zeta^{\prime}_{H}}(s,k+1)}{C_{\zeta^{\prime}_{H}}(s,k)}
=
\frac{2^{k+s}B(s,k)}{2^{k+s+1}B(s,k+1)} = \frac{1}{2} \frac{k+s}{k} 
=  
\frac{1}{2} +\frac{s}{2k} < 1
\]
for $k>s$,
point \ref{control-number-summands-f}) of Definition \ref{set-def} holds 
with $k_{\zeta^{\prime}_{H}}(s) = \max(\lceil s \rceil +1; 3) $.
Differentiating \eqref{difference-Hurwitz}, we obtain 
\begin{equation}
\label{difference-Hurwitz-prime}
\zeta^{\prime}(s,z+1) = \zeta^{\prime}(s,z) + \frac{\log z}{z^{s}},
\end{equation}
for every $s>1$ and $z>0$ and hence point \ref{recursive-f}) of Definition \ref{set-def} holds with 
$g_{\zeta^{\prime}_{H}}(s,z) =   (\log z) z^{-s}$, $s>1$, $z>0$. 
Summarising, $\zeta^\prime(s,x)$, $s>1$ being fixed, $x>0$, is in $\scrF$.
\end{Proof}

It is easy to obtain that $\zeta^{\prime}(s,1/2) = 2^s (\log 2) \zeta(s) + (2^s -1) \zeta^\prime(s)$.
Moreover, we have the following 
\begin{Proposition}[The number of summands for $\zeta^\prime(s,x)$]
Let $s>1$ be fixed. Let further  $x\in (0,1)$, $x\ne 1/2$,
$\acc\in \N$, $\acc\ge 2$,
\begin{equation}
\label{r-Hurwitz-prime-estim}
r_{\zeta^{\prime}_{H}}(s,z,\acc) = 
\max
\Bigl\{
\Bigl\lceil \frac{(\acc+1 + 2 \lceil s\rceil) \log 2 
+ \vert \log ( 1- \vert 1-z\vert ) \vert 
+ \vert \log B(s, \lceil s\rceil +1)\vert}{ \vert \log  \vert 1-z \vert \vert}\Bigr\rceil -1;
\lceil s\rceil+1; 3
\Bigr\},
\end{equation}
and 
$r^\prime_{\zeta^{\prime}_{H}}(s,z,\acc) = r_{\zeta^{\prime}_{H}}(s,1+z,\acc)$.
For $x\in (1/2,1)$ there exists $\theta_1=\theta_1(s,x)\in (-1/2,1/2)$ such that
\begin{equation} 
\label{zeta-Hurwitz-prime-x>1/2}
\zeta^{\prime} (s,x)  
=  
\zeta^{\prime} (s) - \frac{\log x}{x^s} 
+
\sum_{k=1}^{r_{\zeta^{\prime}_{H}}(s,x,\acc)} c_{\zeta^{\prime}_{H}}(s,k) (1-x)^{k}
+\vert \theta_1 \vert 2^{-\acc}.
\end{equation}
For $x\in (0,1/2)$ we have that
there exists $\theta_2=\theta_2(s,x)\in (-1/2,1/2)$ such that
\begin{equation} 
\label{zeta-Hurwitz-prime-x<1/2}
\zeta^{\prime}  (s,x)  
= 
\zeta^{\prime} (s) 
- \frac{\log x}{x^{s}} 
- \frac{\log(1+x)}{(1+x)^s} 
+ \sum_{k=1}^{r^\prime_{\zeta^{\prime}_{H}}(s,x,\acc)} c_{\zeta^{\prime}_{H}}(s,k)  (-x)^{k} 
+ \vert \theta_2 \vert 2^{-\acc}.
\end{equation}
\end{Proposition}
\begin{Proof}
For $r\ge \max(s;2)$ and $s>1$, using \eqref{alpha-estim} and
 Lemma \ref{stima-errore-criterio-rapporto} we have
\begin{equation*}
\Bigl\vert \sum_{k=r+1}^{\infty} 
c_{\zeta^{\prime}_{H}}(s,k)
( 1-z )^{k}
\Bigr\vert
<
1.64
\sum_{k=r+1}^{\infty} 
\frac{\vert 1-z\vert ^k}{2^{k+s}B(s,k)}
<
\frac{0.82}{2^{s+r}B(s,r+1)}
\frac{\vert 1-z \vert^{r+1}}{1-\vert 1-z\vert}.
\end{equation*}
 Let now $\acc\in \N$, $\acc\ge 2$. 
Arguing as in Section \ref{shifting-trick-f},
for $z\in (0,2)$ we obtain that
$\vert \sum_{k=r+1}^{\infty} 
c_{\zeta^{\prime}_{H}}(s,k) 
( 1-z )^{k}
\vert < 2^{-\acc-1}$
for
$r\ge r_{\zeta^{\prime}_{H}}(s,z,\acc)$, as defined in \eqref{r-Hurwitz-prime-estim}.
For $x\in (1/2,1)$, we can compute $\zeta^\prime(s,x)$
by specialising \eqref{f-x>1/2} thus  obtaining \eqref{zeta-Hurwitz-prime-x>1/2}.
For $x\in (0,1/2)$, thanks to \eqref{difference-Hurwitz-prime} we have $g_{\zeta^\prime_H}(s,x) =  x^{-s}(\log x)$
and $r^\prime_{\zeta^\prime_H}(x,\acc) = r_{\zeta^\prime_H}(1+x,\acc)$;
hence  \eqref{f-x<1/2}  becomes \eqref{zeta-Hurwitz-prime-x<1/2}.
\end{Proof} 
   
The reflection formulae are collected in the following
\begin{Proposition}[The reflection formulae for $\zeta^\prime(s,x)$]
\label{DIF-formulae-Hurwitz-prime}
Let $s>1$ be fixed. Let further  $x\in (0,1)$, $x\ne 1/2$,
$\acc\in \N$, $\acc\ge 2$, $r_1(s,x,\acc)=r^\prime_{\zeta^{\prime}_{H}}(s,x,\acc)/2$
and $r_2(s,x,\acc)  =r_{\zeta^{\prime}_{H}}(s,x,\acc)/2$.
Recalling the definition of $c_{\zeta^{\prime}_{H}}(k ,s)$ in \eqref{coeff-hurwitz-prime-def},
there exists $\theta=\theta(s,x)\in (-1/2,1/2)$  such that 
for $0<x <1/2$ we have
\begin{align}
\notag 
\zeta^{\prime}(s,x) +\zeta^{\prime}(s,1-x)
&=
2 \zeta^{\prime}(s)
- \frac{\log x}{x^s}
- \frac{\log (1-x)}{(1-x)^s}
- \frac{\log (1+x)}{(1+x)^s}
+ 
2\sum_{\ell=1}^{r_1}   c_{\zeta^{\prime}_{H}}(s,2\ell) x^{2\ell}  
+ \vert \theta \vert 2^{-\acc},
\\
\zeta^{\prime}(s,x) -\zeta^{\prime}(s,1-x)
&=
- \frac{\log x}{x^s}
+ \frac{\log (1-x)}{(1-x)^s}
- \frac{\log (1+x)}{(1+x)^s}
\label{zeta-Hurwitz-prime-DIF-odd-x<1/2}
-
2\sum_{\ell=1}^{r_1}
c_{\zeta^{\prime}_{H}}(s,2\ell-1)  x^{2\ell-1} 
+ \vert \theta \vert 2^{-\acc},
\end{align}
and for $1/2<x <1$ we have
\begin{align*}
\notag
\zeta^{\prime}(s,x) +\zeta^{\prime}(s,1-x)
&
=
2 \zeta^{\prime}(s)
- \frac{\log x}{x^s}
- \frac{\log (1-x)}{(1-x)^s}
- \frac{\log (2-x)}{(2-x)^s} 
+
2\sum_{\ell=1}^{r_2}   c_{\zeta^{\prime}_{H}}(s,2\ell) (1-x)^{2\ell} 
+ \vert \theta \vert 2^{-\acc} ,
\\ 
\notag
\zeta^{\prime}(s,x) -\zeta^{\prime}(s,1-x)
&=
-\frac{\log x}{x^s}
+ \frac{\log (1-x)}{(1-x)^s}
+ \frac{\log (2-x)}{(2-x)^s}
+
2\sum_{\ell=1}^{r_2}  c_{\zeta^{\prime}_{H}}(s,2\ell-1) (1-x)^{2\ell-1}  
+ \vert \theta \vert 2^{-\acc} .
\end{align*}
\end{Proposition}
Similar formulae hold for the infinite series too.

For this function we can apply Lemma \ref{Euler-Maclaurin} to study its tail but 
the final formulae are much more complicated than in the previous cases 
due to the form of the derivatives of $g_{\zeta^\prime_H}(s,x) = x^{-s}(\log x)$.
Since in  Section \ref{DirichletL} we will just need the values 
of $\zeta^{\prime}(s,x)$ for $x\in (0,1)$, we do not insert this topic here.
\begin{Remark} [Coefficients computations] 
In this case we have to precompute both $\zeta(k+s)$ and  $\zeta^{\prime}(k+s)$;
moreover
the  coefficients
$b_{s,k} :=  
(k B(s,k))^{-1}$
can be computed as in \eqref{recursive-coeff-hurwitz}.
For $d_{s,k} := \psi(k+s)- \psi(s)$, using \eqref{difference-psi} we obtain
$d_{s,0} = 0$, $d_{s,1} = s^{-1}$,  $d_{s,2} = (s+1)^{-1}+d_{s,1}$, 
and hence 
\begin{equation*}
d_{s,k+1} = \psi(k+s+1)- \psi(s)  
= \frac{1}{s+k} + \psi(k+s)- \psi(s) 
= \frac{1}{s+k} +d_{s,k} .
\end{equation*}
\end{Remark}

We finally remark that  the maximal order of magnitude for $\zeta^{\prime}(s,x)$
is $x^{-s}(\log x)$ for $x\to 0^+$; so to get  its accurate evaluation 
one needs at least  $\lceil s \vert \log_2 x \vert \rceil 
+ \vert \log_2 \vert \log x \vert\vert \rceil$ binary digits.

\subsection{The Dirichlet $\beta$-function}
\label{Dirichlet-beta-def-sec}
The results in Proposition \ref{DIF-formulae-Hurwitz} 
reveal that we have an efficient strategy 
to compute the cases in which are involved the  Hurwitz zeta-function values 
at the $(s,x)$ and $(s,1-x)$-points. 
In fact there exists a famous example of a function of this kind,
the Dirichlet $\beta$-function:
\[
\beta(s) := 
\frac{\zeta(s,1/4) - \zeta(s,3/4)}{4^s}
=
\sum_{n=0}^{\infty} \frac{(-1)^n}{(2n+1)^{s}}
=
L(s,\chi_{-4}),
\] 
where $\chi_{-4}$ is 
the quadratic Dirichlet character mod $4$.
Recalling that the Euler beta-function is denoted as $B(x,y)$,
from  \eqref{enlarged-radius-Hurwitz} and \eqref{difference-Hurwitz}, 
for $x\in (0,1/2)$ and $s>1$ we have
\begin{equation*}
\zeta(s,x) -\zeta(s,1-x)
=
\frac{1}{x^s} +  \frac{1}{(1+x)^s} -  \frac{1}{(1-x)^s}
-
2\sum_{\ell=1}^{\infty}
\frac{\zeta(2\ell-1+s) -1}{(2\ell-1) B(s,2\ell-1)}  x^{2\ell-1},
\end{equation*}
which is the infinite series version of \eqref{zeta-Hurwitz-DIF-odd-x<1/2}.
Hence, for $s>1$ we  obtain that
\begin{equation}
\label{dirichlet-beta-def}
\beta(s) 
= 
1
-  \frac{1}{3^s} + \frac{1}{5^s} 
-
\frac{1}{2^{2s-3}}   \sum_{\ell=1}^{\infty}
\frac{\zeta(2\ell-1+s)-1}{ 16^{\ell}(2\ell-1) B(s,2\ell-1)} 
\end{equation}
and a corresponding truncated formula can be obtained arguing
as we did in Proposition \ref{DIF-formulae-Hurwitz} to obtain \eqref{zeta-Hurwitz-DIF-odd-x<1/2}.
Recalling \eqref{coeff-Hurwitz-def} and \eqref{coeff-hurwitz-prime-def}, or 
differentiating both sides of formula \eqref{dirichlet-beta-def}, for $s>1$ we can also get 
\begin{align*}
\beta^\prime(s) 
&= 
-2 (\log 2) \beta(s)
+
4^{-s} \bigl(
\zeta^{\prime}(s,1/4) -\zeta^{\prime}(s,3/4)\bigr)
\\
&=
-2 \log 2
\Bigl(1
+ \frac{1}{5^s} -  \frac{1}{3^s}
-
\frac{1}{2^{2s-3}}  \sum_{\ell=1}^{\infty}
c_{\zeta_{H}}(s, 2\ell-1) 
16^{-\ell} \Bigr)
\\& \hskip1cm
+
2 \log 2
- \frac{\log(5/4)}{5^s}
+\frac{\log(3/4)}{3^s}
-
\frac{1}{2^{2s-3}}  \sum_{\ell=1}^{\infty}
c_{\zeta^{\prime}_{H}}(s, 2\ell-1) 
16^{-\ell}
\\&
=
\frac{\log 3}{3^s}
-\frac{\log 5}{5^s} 
- \frac{1}{2^{2s-3}} \sum_{\ell=1}^{\infty}
\frac{ (\zeta(2\ell-1+s)-1) (\psi(2\ell-1+s)
- 
\psi(s) - 2 \log 2) + \zeta^{\prime}(2\ell-1+s)}
{16^{\ell}(2\ell-1) B(s,2\ell-1)}
\end{align*}
and a corresponding formula involving a truncated sum can be obtained arguing
as we did in Proposition \ref{DIF-formulae-Hurwitz-prime}.
Moreover, using \eqref{zeta-Hurwitz-prime-DIF-odd-x<1/2} and \eqref{dirichlet-beta-def}, the logarithmic derivative of $\beta(s)$
can be expressed as follows:
\begin{align*}
\frac{\beta^\prime(s) }{\beta(s) } &= -2 \log2 +  \frac{\zeta^{\prime}(s,1/4) -\zeta^{\prime}(s,3/4)}{\zeta(s,1/4) -\zeta(s,3/4)}
\\
&= -2 \log2  +  
\frac
{2\log 2 \bigl(1-3^{-s}+5^{-s} \bigr) +(\log 3)3^{-s} - (\log 5)5^{-s} -2^{3-2s}\sum_{\ell=1}^{\infty} c_{\zeta^{\prime}_{H}}(s, 2\ell-1) 16^{-\ell} }
{1-3^{-s}+5^{-s} -2^{3-2s}\sum_{\ell=1}^{\infty} c_{\zeta_{H}}(s, 2\ell-1) 16^{-\ell}},
\end{align*}
where $c_{\zeta_{H}}(s,k)$ and $c_{\zeta^{\prime}_{H}}(s,k)$ are respectively defined in
\eqref{coeff-Hurwitz-def} and \eqref{coeff-hurwitz-prime-def}.

In Section \ref{tests} we will show that this approach gives pretty good performances
when compared with a standard implementation of $\beta(s)$, $\beta^\prime(s)$
and $\beta^\prime(s)/\beta(s)$.

\subsection{The Catalan constant}
Another famous quantity connected with the values of the Hurwitz
zeta-function is the Catalan constant. 
Specialising what we wrote for the Dirichlet $\beta$-function, we
will obtain a  new formula for $G$ which is very similar to 
Glaisher's result \cite{Glaisher1912} but has a better convergence
speed. In fact, there exist much faster \emph{ad hoc} algorithms to compute
$G$ and so we insert this paragraph here just to show how our
way of arguing can lead to obtain new formulas for classical quantities.
Using \eqref{dirichlet-beta-def}  for  $s=2$, we  obtain that
\begin{equation}
\label{similar-Glaisher-catalan}
G 
: =  
\beta(2) 
= 
\frac{209}{225}
- \frac{1}{2} \sum_{\ell=1}^{\infty} 
\frac{\zeta(2\ell+1)-1}{16^{\ell}(2\ell-1) B(2,2\ell-1)} 
=
\frac{209}{225} - \sum_{\ell=1}^{\infty} 
\frac{\ell}{16^{\ell}} (\zeta(2\ell+1) -1)
\end{equation}
in which we also used 
$d_{k,2} =  (kB(2,k))^{-1}= k+1$, for $k\ge 1$,
and $B(x,y)$ is the Euler beta-function.
Moreover, thanks to Lemma \ref{elementary-estim}, the order of magnitude
of the summands is, roughly speaking,
about $\ell \cdot 64^{-\ell}$; 
much smaller than in Glaisher's formula.
We also remark that
a truncated formula for \eqref{similar-Glaisher-catalan} can be proved directly 
or using \eqref{zeta-Hurwitz-DIF-odd-x<1/2}.
In a multiprecision computation of $G$, the bottleneck in using \eqref{similar-Glaisher-catalan} is 
the large number of high-precision Riemann zeta-function values 
needed, since in this case they are just used once.
We finally recall that Glaisher
was able to obtain $G$ with an accuracy of $32$ decimal digits;
currently $G$ is now known up to $6\cdot10^{11}$ decimal digits using much faster
\emph{ad hoc} algorithms, see, e.g., Kim \cite{Kim2019}.

\section{How to compute the Dirichlet $L$-functions for $s\ge 1$}  
\label{DirichletL}
Assume that  $q$ is an odd prime and let $\chi$ be a non-principal Dirichlet character mod $q$.
For $s=1$ the approach is the one already published in \cite{Languasco2021},
\cite{Languasco2021a} and \cite{LanguascoR2021}. 
For $s>1$ will use the well known formulae  
\begin{equation}
\label{L-hurwitz} 
L(s, \chi) = q^{-s} \sum_{a=1}^{q-1} \chi(a)\zeta\Bigl(s,\frac{a}{q}\Bigr)
\quad\textrm{and}
\quad
L^{\prime}(s, \chi) = -(\log q) L(s,\chi) 
+ q^{-s} \sum_{a=1}^{q-1} \chi(a)\zeta^{\prime}\Bigl(s,\frac{a}{q}\Bigr);
\end{equation}
for the first one see, e.g., Cohen \cite[Proposition 10.2.5]{Cohen2007}, while the second one
can be immediately obtained by differentiating the first.
It is clear that we need the values of $\zeta(s, z)$ and  $\zeta^{\prime}(s, z)$ 
for $z =a/q\in (0,1)$. To this goal we can respectively use the results described in Sections 
\ref{Hurwitz-zeta-section} and \ref{Hurwitz-prime-section};
in this application it is particularly efficient the fact that, 
for $s>1$ fixed, both the coefficients $c_{\zeta_{H}}(s,k)$ in \eqref{coeff-Hurwitz-def} and 
$c_{\zeta^{\prime}_{H}}(s,k)$ in \eqref{coeff-hurwitz-prime-def} 
can be precomputed.  

If $q$ is not small, in fact, the computational cost of performing the sums over 
$a$ in \eqref{L-hurwitz} becomes too large. But the trivial
summation procedure can be replaced by the use of the Fast Fourier Transform (FFT)
algorithm; it is in this context that the reflection formulae of Propositions
\ref{DIF-formulae-Hurwitz} and \ref{DIF-formulae-Hurwitz-prime}
will be useful. 

We recall that the  FFT-procedure  is a quite fast, but memory demanding, algorithm. 
It computes a linear combination
of complex exponentials whose coefficients
are the values of a given finite sequence $\scrA$ of $N$ complex numbers.
Instead of performing such a summation term by term,
which would lead to a total computational cost of $\Odi{N^2}$ products,
the FFT procedure implements a  \emph{divide et impera} strategy 
that uses the \emph{decimation in time} or the
\emph{decimation in frequency} ideas, see Sections  \ref{DIT-sect}-\ref{DIF-sect}.
This reduces the total computational cost to $\Odi{N \log N}$ products
but requires the storage of at least one copy of the whole sequence $\scrA$.
Hence, roughly speaking, we can say that $N$
memory positions are required to perform such a computation;
in practice, more memory space is in fact needed to keep
track of the several steps  an implementation of the FFT requires.
We also recall that it is not an easy task to efficiently and accurately 
implement the FFT algorithm; we refer to the papers of Cooley-Tukey \cite{CooleyT1965},
Cochran, Cooley et al.~\cite{CochranC1967},
Rader \cite{Rader1968} and to Arndt's book \cite[Part III]{Arndt2011} for more details.
We show now how to use the FFT to compute the values of 
$L(s, \chi)$ and $L^\prime(s, \chi)$.

\subsection{The Fast Fourier Transform setting}
\label{FFT-setting}
First of all we have to introduce the following
\begin{Definition}[The Discrete Fourier Transform (DFT)]
\label{FFT-sum-def}
Let $N\in \N$, $\scrA$ be a sequence whose elements $\scrA_k \in \C$, $k=0,\dotsc, N-1$. We define
the \emph{Discrete Fourier Transform} of $\scrA$ as the  sequence $\boldF_{\sigma}(\scrA)$ such that
\begin{equation}
\label{DFT-sum}
 \boldF_{\sigma}(\scrA)_j : =
\sum_{k=0}^{N-1}  e\Bigl(\frac{\sigma  j k}{q-1}\Bigr) 
\scrA_k,
\end{equation}
where $e(x):=\exp(2\pi i x)$, $j\in\{0,\dotsc,N-1\}$, and $\sigma=\pm 1$ is a fixed parameter.
\end{Definition}

The sequence $\boldF_{1}(\scrA)$ is usually called  the \emph{forward} DFT
of $\scrA$, while, in the other case, it is called the \emph{backward} 
(or \emph{inverse}) DFT of $\scrA$. 
We want now to connect the summations in \eqref{L-hurwitz}
with $\boldF_{\sigma}(\scrA)$. The first ingredient we need is the following 
lemma whose proof can be found, e.g., in Davenport \cite[Chapter 4]{Davenport2000}.
\begin{Lemma}[Dirichlet characters representation]
\label{Dirichlet-character-rep}
Let $q$ be an odd prime and let $g$ be a primitive root of $q$.
Let further $\chi_1$ be the Dirichlet character mod $q$ given by
$\chi_1(g) = e^{2\pi i/(q-1)}$. Then the set of the Dirichlet characters
mod $q$ is $X:=\{\chi_1^j \colon j=0,\dotsc,q-2\}$. Moreover, 
the principal character, $\chi_0 \bmod q$, corresponds to $j=0$ 
in the previous representation.
\end{Lemma} 
It is well known that the problem of finding a primitive root $g$ of $q$ is a 
computationally hard one, see, e.g., Shoup \cite[Chapters 11.1-11.4]{Shoup2005},
but, for each fixed prime $q$, we need to find it just once. In the applications,
for each involved prime $q$ we can save such a $g$ and reuse
it every time we have to work again mod $q$.

For every $k\in \{0,\dotsc,q-2\}$, let now denote $g^k\equiv a_k\in\{1,\dotsc,q-1\}$.
We will use \eqref{DFT-sum}, with $N=q-1$, in the cases in which
$\scrA_k=f(a_k/q)$, and $f(\cdot)$ is either $\zeta (s,\cdot)$ or $\zeta^\prime (s,\cdot)$.
We also remark that, for a fixed odd prime $q$, we do not need to generate the sequence $a_k/q$ 
several times (one for  $\zeta (s,a_k/q)$ and one for  $\zeta^\prime(s,a_k/q)$, for example), but just once.

\subsection{FFT: decimation in time (DIT)}
\label{DIT-sect}
Now we recall one of the two main  ideas used in the Cooley-Tukey \cite{CooleyT1965} FFT algorithm:
the \emph{decimation in time} strategy. We describe only the special case in which the
length of the transform is even;  
in fact the FFT algorithm can be used in the general case too,
we refer to  \cite{CooleyT1965}, \cite{CochranC1967}, \cite{Rader1968} and to \cite[Part III]{Arndt2011} for more details.

In the following we will use $\scrE$ and $\scrO$ to respectively denote the subsequences $(\scrA)_{2k}$ and  
$(\scrA)_{2k+1}$ of $\scrA$.
We show now that the first $\overline{q}=(q-1)/2$ elements 
of the sequence $\boldF_{\sigma}(\scrA)$ (also called  \emph{the left part of} $\boldF_{\sigma}(\scrA)$) 
can be written using
the sequences $\boldF_{\sigma}(\scrE)$ and $\boldF_{\sigma}(\scrO)$; the same also holds
for the second $\overline{q}$ elements  of $\boldF_{\sigma}(\scrA)$ (also called  \emph{the right part of} $\boldF_{\sigma}(\scrA)$).

\begin{Lemma}[FFT: decimation in time]
\label{DIT-lemma}
Let $q$ be an odd prime and $\overline{q}=(q-1)/2$. Let $k \in\{0,\dotsc,q-2\}$ and $\scrA$ be a 
sequence having $q-1$ elements $\scrA_k \in \C$. Let  further
$\scrE, \scrO$ respectively denote the subsequences $(\scrA)_{2k}$,  $(\scrA)_{2k+1}$ of $\scrA$ and
 $\boldF_{\sigma}(\scrA)$ be defined as in \eqref{DFT-sum}
where $e(x) =\exp(2\pi i x)$ and $\sigma=\pm 1$.
Then for every $j\in\{0,\dotsc,\overline{q}-1\}$ we have
\begin{equation}
\label{DIT-left}
 \boldF_{\sigma}(\scrA)_j 
=
( \boldF_{\sigma}(\scrE))_j  +  e\Bigl(\frac{\sigma j}{q-1}\Bigr)    ( \boldF_{\sigma}(\scrO))_j.
\end{equation}
Moreover, for every $j\in\{\overline{q}, q-2\}$
we have
\begin{equation}
\label{DIT-right}
 \boldF_{\sigma}(\scrA)_{j} 
=
( \boldF_{\sigma}(\scrE))_{j-\overline{q}}  -  e\Bigl(\frac{\sigma (j-\overline{q})}{q-1}\Bigr)    ( \boldF_{\sigma}(\scrO))_{j-\overline{q}}.
\end{equation}
\end{Lemma}
\begin{Proof}
Let $j=t+\delta \overline{q}$, where $\delta \in\{0,1\}$. Taking $k=0,\dotsc, \overline{q}-1$, 
and splitting the sum in \eqref{DFT-sum} according to parity we obtain 
\begin{equation}
\label{DIT-splitting}
\boldF_{\sigma}(\scrA)_j 
 =
\sum_{k=0}^{\overline{q}-1}  
e\Bigl(\frac{\sigma (t+\delta \overline{q}) 2k}{q-1}\Bigr) \scrA_{2k}
+
e\Bigl(\frac{\sigma (t+\delta \overline{q})}{q-1}\Bigr)  
\sum_{k=0}^{\overline{q}-1}   
e\Bigl(\frac{\sigma (t+\delta \overline{q}) 2k}{q-1}\Bigr) \scrA_{2k+1}. 
\end{equation}

The left part of $\boldF_{\sigma}(\scrA)$ has $\delta = 0$ in \eqref{DIT-splitting}. In this case 
$j=t\in\{0,\dotsc,\overline{q}-1\}$ and
since $2\overline{q} = q-1$,  we have 
$e\bigl(\frac{\sigma j 2k}{q-1}\bigr) = e\bigl(\frac{\sigma j k}{\overline{q}}\bigr)$
for every $k\in \N$.
Hence \eqref{DIT-splitting} becomes 
\[
\boldF_{\sigma}(\scrA)_j 
=
\sum_{k=0}^{\overline{q}-1}  e\Bigl(\frac{\sigma j k}{\overline{q}}\Bigr) \scrA_{2k}
+
e\Bigl(\frac{\sigma j}{q-1}\Bigr)  
\sum_{k=0}^{\overline{q}-1}   e\Bigl(\frac{\sigma j k}{\overline{q}}\Bigr) \scrA_{2k+1}.
\]
Clearly both the sums are of the type in  \eqref{DFT-sum} but their lengths are now 
$\overline{q}$ instead of $q-1$;  the input sequence for the first one is $\scrE$ 
while for the second is $\scrO$. Hence \eqref{DIT-left} follows. 

The right part of $\boldF_{\sigma}(\scrA)$ has  $\delta = 1$ in \eqref{DIT-splitting}. 
In this case $j=t+\overline{q}\in\{\overline{q},\dotsc,q-2\}$ and
since $2\overline{q} = q-1$,  we have 
$e\bigl(\frac{\sigma (t+ \overline{q}) 2k}{q-1}\bigr)
= 
e\bigl(\frac{\sigma t (2k)}{q-1}\bigr)  e(\sigma k) 
=
e\bigl(\frac{\sigma t k}{\overline{q}}\bigr)$ 
for every $k\in \N$. Moreover, we also have 
$e\bigl(\frac{\sigma (t+ \overline{q})}{q-1}\bigr) 
= 
e\bigl(\frac{\sigma t}{q-1}\bigr) e\bigl(\frac{\sigma}{2}\bigr)  
= 
- e\bigl(\frac{\sigma t}{q-1}\bigr)$.
Hence \eqref{DIT-splitting} becomes 
\[
 \boldF_{\sigma}(\scrA)_{t+\overline{q}} 
=
\sum_{k=0}^{\overline{q}-1}  
e\Bigl(\frac{\sigma t k}{\overline{q}}\Bigr) \scrA_{2k}
-
e\Bigl(\frac{\sigma t}{q-1}\Bigr)  
\sum_{k=0}^{\overline{q}-1}   
e\Bigl(\frac{\sigma t k}{\overline{q}}\Bigr) \scrA_{2k+1}.
\]
Clearly both the sums in the previous equation are of the type in  \eqref{DFT-sum} but their lengths are now 
$\overline{q}$ instead of $q-1$; the input sequence for the first one is $\scrE$ 
while for the second is $\scrO$. Hence, recalling $j=t+\overline{q}$, \eqref{DIT-right} follows.
\end{Proof}

Lemma \ref{DIT-lemma} shows that both the left and the right parts of  $\boldF_{\sigma}(\scrA)$
are suitable combinations of $\boldF_{\sigma}(\scrE)$ and $\boldF_{\sigma}(\scrO)$.
Since  both $\scrE$ and $\scrO$ have half a length of the original sequence $\scrA$,
we have reduced the problem of computing a transform of a sequence having $q-1$ elements
to the problem of computing two transforms of sequences having  $(q-1)/2$ elements each.
This, together with what we will see in the next section, 
is the starting point of a recursive procedure, called the \emph{Fast Fourier Transform}
algorithm, that leads to compute 
$\boldF_{\sigma}(\scrA)$, a transform of length $q-1$, in 
$\Odi{\log q}$ steps and $\Odi{q\log q}$ products and sums.
We refer to Cochran, Cooley et al.~\cite{CochranC1967}, and to Arndt's book \cite[Part III]{Arndt2011}, 
for more on this topic.

We can use the decimation in time strategy, but 
the parity of the Dirichlet characters, see Lemma \ref{Dirichlet-character-rep}, 
can be detected by working on the parity of the indices of $\boldF_{\sigma}(\scrA)$, not of $\scrA$.
Luckily, another strategy is possible; we will show it in the next section.

\subsection{FFT: decimation in frequency (DIF)}
\label{DIF-sect}
We can now show how the \emph{decimation in frequency} 
strategy works.  We assume that 
in \eqref{DFT-sum} one has to distinguish between the parity of $j$
(hence on the parity of the Dirichlet characters, see Lemma \ref{Dirichlet-character-rep}).
In the next lemma we will obtain that the subsequences
of $\boldF_{\sigma}(\scrA)$  having odd and even indices can be respectively obtained with a DFT
applied onto some suitable modifications of the left and right parts of $\scrA$.

\begin{Lemma}[FFT: decimation in frequency]
\label{DIF-Lemma}
Let $q$ be an odd prime and $\overline{q}=(q-1)/2$. Let $k \in \{0,\dotsc, q-2\}$ and $\scrA$  be a 
sequence having $q-1$ elements  $\scrA_k \in \C$. Let   $e(x) =\exp(2\pi i x)$,
$\sigma=\pm 1$ and
$\boldF_{\sigma}(\scrA)$  be defined as in \eqref{DFT-sum}. For
$k \in \{0,\dotsc,  \overline{q}-1\}$ define the sequences $\scrB$ and $\scrC$
whose elements are
\begin{equation}
\label{Bk-Ck-defs}
\scrB_k :=
\scrA_{k} + \scrA_{k+\overline{q}} 
\quad
\text{and}
\quad
\scrC_k := 
e\Bigl(\frac{\sigma k}{q-1}\Bigr)   
\bigl(  \scrA_{k} - \scrA_{k+\overline{q}}   \bigr).
\end{equation}
%
Let further $t\in \{0,\dotsc, \overline{q}-1\}$. 
We have
\begin{equation}
\label{DIF-final}
 \boldF_{\sigma}(\scrA)_{2t}
 =  \boldF_{\sigma}(\scrB)_t
 \quad
 \textrm{and}
 \quad
  \boldF_{\sigma}(\scrA)_{2t+1}
 =  \boldF_{\sigma}(\scrC)_t.
\end{equation}

\end{Lemma}
\begin{Proof}
Recalling  $\sigma=\pm 1$ and $\overline{q} =(q-1)/2$, by splitting in two halves the sum over $k$  
in \eqref{DFT-sum}, for every $j=0,\dotsc,q-2$  we have that
\begin{equation*}
 \boldF_{\sigma}(\scrA)_j 
= 
\sum_{k=0}^{\overline{q}-1} 
e\Bigl(\frac{\sigma  j k}{q-1}\Bigr)  
\scrA_k 
+
\sum_{k=0}^{\overline{q}-1} 
e\Bigl(\frac{\sigma  j (k+\overline{q})}{q-1}\Bigr)
\scrA_{k+\overline{q}}
=
\sum_{k=0}^{\overline{q}-1}
e\Bigl(\frac{\sigma  j k}{q-1}\Bigr)
\bigl(
\scrA_k 
+
(-1)^{j} 
\scrA_{k+\overline{q}}
\bigr),
\end{equation*}
in which we used  that
$e\bigl(\frac{\sigma  j (k+\overline{q})}{q-1}\bigr)
=
e\bigl(\frac{\sigma  j k}{q-1}\bigr)
e\bigl(\frac{\sigma  j }{2}\bigr) 
=
(-1)^j
e\bigl(\frac{\sigma  j k}{q-1}\bigr)$.
Let now $j=2t+\ell$, where $\ell\in\{0,1\}$ and $t=0,\dotsc, \overline{q}-1$. 
Then the previous equation becomes
\begin{equation}
\label{DIF}
 \boldF_{\sigma}(\scrA)_{2t+\ell}
=
\sum_{k=0}^{\overline{q}-1}
e\Bigl(\frac{\sigma  t k}{\overline{q}}\Bigr)   
e\Bigl(\frac{\sigma \ell k}{q-1}\Bigr)   
\bigl( 
\scrA_k 
+
(-1)^{\ell} 
\scrA_{k+\overline{q}}
\bigr) 
=
\begin{cases}
\sum\limits_{k=0}^{\overline{q}-1} 
e\bigl(\frac{\sigma t k}{\overline{q}}\bigr) \scrB_k 
& \textrm{if} \ \ell =0\\
\sum\limits_{k=0}^{\overline{q}-1} 
e\bigl(\frac{\sigma t k}{\overline{q}}\bigr) \scrC_k  
& \textrm{if} \ \ell =1,\\
\end{cases}
\end{equation}
where $\scrB_k$ and $\scrC_k$ are defined as in \eqref{Bk-Ck-defs}.
Clearly both the sums on the right hand side of \eqref{DIF} are of the type in  \eqref{DFT-sum} 
but their lengths are now 
$\overline{q}$ instead of $q-1$; the input sequence for the first one is $\scrB$ 
while for the second one it is  $\scrC$. Hence \eqref{DIF-final} follows.
\end{Proof}

Remark that both $\scrB$ and $\scrC$ in \eqref{Bk-Ck-defs} are built using the 
left and the right parts of $\scrA$.

Recalling that the set of the Dirichlet characters are represented as in Lemma \ref{Dirichlet-character-rep},
Lemma \ref{DIF-Lemma} hence shows how to  split the 
original problem according to their parity;
in this way instead of computing a DFT of length $q-1$
we can evaluate two DFTs of length $(q-1)/2$ each, applied on suitably
modified sequences according to \eqref{DIF-final}.  
%
We show now how to insert the reflection formulae for  $f(\cdot)$ 
in this setting.

\begin{Lemma}[FFT-DIF and the reflection formulae for $f$]
\label{DIF-reflection-Lemma}
Let $q$ be an odd prime, $g$ be a primitive root of $q$
and $f\colon (0,1) \to \C$ be a function.
Let  $\scrA_k= f(a_k/q)$,
where $a_k \equiv g^k \bmod q \in\{1,\dotsc, q-1\}$,  $k=0,\dotsc,q-2$,
and $\scrB$, $\scrC$ be defined as in \eqref{Bk-Ck-defs}.
Let further $\overline{q}=(q-1)/2$.
%
Then, for every $k \in \{0,\dotsc,  \overline{q}-1\}$, we have
\begin{equation}
\label{DIF-reflection}
\scrB_k :=
f\Bigl(\frac{a_k}{q}\Bigr) 
+  
f\Bigl(1-\frac{a_{k}}{q}\Bigr) 
\quad
\textrm{and}
\quad
\scrC_k := 
e\Bigl(\frac{\sigma k}{q-1}\Bigr)   
\Bigl(  f\Bigl(\frac{a_k}{q}\Bigr) 
- 
f\Bigl(1-\frac{a_{k}}{q}\Bigr)  \Bigr).
\end{equation}
\end{Lemma}
\begin{Proof}
Since $g$ is a primitive root of $q$, it trivially follows that 
$g^{\overline{q}} \equiv -1 \bmod{q}$,
where $\overline{q}=(q-1)/2$.
Hence 
we have
\(
a_{k+\overline{q}} \equiv g^{k+\overline{q}}  
\equiv 
- a_k  \equiv q-a_k \bmod{q}
\) 
thus obtaining
\begin{equation*} 
\scrA_{k+\overline{q}}
=
f \Bigl(\frac{a_{k+\overline{q}}}{q}\Bigr)  
= 
f\Bigl(\frac{q-a_{k}}{q}\Bigr)
=
f\Bigl(1-\frac{a_{k}}{q}\Bigr).
\end{equation*}
The lemma is hence proved by recalling the definitions
of $\scrB_k$ and $\scrC_k$ in  \eqref{Bk-Ck-defs}.
\end{Proof}

Equation \eqref{DIF-reflection} of
Lemma \ref{DIF-reflection-Lemma} clearly shows that the sequences involved
in the decimation in frequency procedure of Lemma \ref{DIF-Lemma} can be 
computed using the reflection formulae for $f(\cdot)$.
Recalling that in our application $f(\cdot)$ is either $\zeta (s,\cdot)$ or $\zeta^\prime (s,\cdot)$,
Propositions \ref{DIF-formulae-Hurwitz} and \ref{DIF-formulae-Hurwitz-prime} provide what we need
to handle the computation of $L(s,\chi)$, $s>1$, using Lemmas \ref{DIF-Lemma}-\ref{DIF-reflection-Lemma},
or, in other words, using the FFT-DIF procedure.
In other cases, for example for computing the values of the Dirichlet $L$-function and its first derivative
at $s=1$, we can apply Lemmas \ref{DIF-Lemma}-\ref{DIF-reflection-Lemma} to the function  $f(\cdot)$ 
needed in those cases; this is  what we did in \cite{Languasco2021},
\cite{Languasco2021a} and \cite{LanguascoR2021}. 

\subsection{The gain in the computational effort due to the use of the reflection formulae in the decimation in frequency strategy}
\label{further-FFT-gain}
Let $\Delta\in(0,1)$ be the accuracy we would like to achieve in computing the input sequence $\scrA$
of the FFT procedure.
Recall that $q$ is an odd prime, $g$ is a primitive root mod $q$ and  $a_k \equiv g^k \bmod q$, $a_k\in\{1,\dotsc,q-1\}$.
For  $\scrA_k = f(a_k/q)$ equal to either $\zeta (s,a_k/q)$ or $\zeta^\prime (s,a_k/q)$, $s>1$ being fixed,
but also for the analogous sequences that can be obtained using $\log\Gamma$ and $\psi$,  
the number of the required summands to obtain a $\Delta$-accuracy computation 
for $\scrB_k$ and $\scrC_k$ is reduced by a factor of $2$
if compared with the analogous number of summands for $\scrA_k$:
this depends on the fact the summands in
Propositions \ref{DIF-formulae-Hurwitz} and \ref{DIF-formulae-Hurwitz-prime} 
run only over odd, or  even, indices.

Moreover, recalling that the length of the sequences $\scrB ,\scrC$ is half the length of $\scrA$,
the whole computational effort required to start the FFT procedure is then reduced by a factor of $4$.
And the number of products required to perform the whole FFT is reduced by a factor
of at least $2$ (compare the order of magnitude of $q \log q$ with the one of $q/2 \log (q/2)$).
We also remark that the memory usage  required to perform the 
FFT algorithm is reduced by a factor of $2$
since it is now enough to work with a sequence ($\scrB$, for example) whose length is half the original one ($\scrA$).
We further remark that, to perform the FFT-transform of $\scrC$, we can reuse the same memory space used for 
$\scrB$. 

Thus, the combination of the gain in speed and in the memory usage let us work with larger values of $q$. Some
examples of this fact are shown in our works \cite{Languasco2021},
\cite{Languasco2021a}, \cite{LanguascoR2021} in which we used these ideas 
to compute the values of  $L^\prime(1,\chi)/L(1,\chi)$
and other related quantities for every odd prime $q\le 10^7$.

\section{Beyond the set  $\scrF$}
\label{beyond-F} 
\subsection{Introduction}
\label{intro-beyond-F}
We discuss here two further examples: the Bateman $G$-function
and the alternating Hurwitz zeta-function (also called 
the  Hurwitz-type Euler zeta-function, or the Hurwitz $\eta$-function). 
Both functions are defined in terms of functions belonging to $\scrF$,
namely
\begin{equation}
\label{Bateman-G-def}
G(z) : = \psi\Bigl(\frac{z+1}{2}\Bigr) - \psi\Bigl(\frac{z}{2}\Bigr)
= 2 \sum_{n=0}^{\infty} \frac{(-1)^n}{n+z},
\end{equation}
for every $z>0$, where the second relation comes 
from eq.~(1.8.6) on page 20 of 
\cite{ErdelyMOT1953},
and
\begin{equation} 
\label{alt-Hurwitz-def}
\eta(s,z) : =  
2^{-s} \Bigl( \zeta\Bigl(s,\frac{z}{2}\Bigr) - \zeta\Bigl(s,\frac{1+z}{2}\Bigr)\Bigr)
=
\sum_{n=0}^{\infty} \frac{(-1)^n}{(n+z)^s},
\end{equation}
for every $s > 1$ and $z>0$, where the first relation  
follows  by splitting the summands in the series according to parity.
We also have functional equations for both $G(z)$ and $\eta(s,z)$, namely 
\begin{equation} 
\label{difference-G-eta}
G(1+z) = \frac{2}{z} - G(z),
\quad z>0, 
\quad \textrm{and} \quad
\eta(s,1+z) = \frac{1}{z^{s}} - \eta(s,z),
\quad z>0, s>1,
\end{equation}
but  they can just be used to build suitable reflection formulae and not
to handle their tails as we did for the functions belonging to $\scrF$
since they connect $f(s,z+1)$ with $-f(s,z)$ and not to $f(s,z)$ as required
in point \ref{recursive-f}) of Definition \ref{set-def}.
In fact, letting $x>0$ such that $\lfloor x \rfloor   \ge 1$ and  $\{ x\} >0$,
the relations in
\eqref{difference-G-eta}
respectively lead to the following analogues of \eqref{main-recursive-f}:
\[
G(x)
=
(-1)^{\lfloor x \rfloor}
\Bigl(
G(\{ x\})
-
2 \sum_{j=0}^{ \lfloor x\rfloor -1} 
\frac{(-1)^j}{\{ x\} + j}
\Bigr)
\quad 
\textrm{and}
\quad
\eta(s,x)
=
(-1)^{\lfloor x \rfloor}
\Bigl(
\eta(s,\{ x\})
-
\sum_{j=0}^{ \lfloor x\rfloor -1} 
\frac{(-1)^j}{(\{ x\} + j)^s}
\Bigr)
\]
that cannot be handled using the Euler-Maclaurin formula of Lemma \ref{Euler-Maclaurin}
since in both cases we have an alternating sum.
To overcome this obstruction, we will use twice the functional equations for $\psi(z)$ and $\zeta(s,z)$ instead,
the first time on $z/2$ and the second one on $(1+z)/2$.
The final outcome, see Propositions \ref{Bateman-recursive} and \ref{eta-recursive}, will connect $f(s,z)$ to  $f(s,2\{z/2\})$
or to  $-f(s,2\{z/2\}-1)$ according to the fact that 
$z\in (2\ell, 2\ell+1)$ or $z\in (2\ell+1, 2\ell+2)$, $\ell \in \N$,
and the tails thus obtained
will be evaluated using  twice the corresponding treatment
for the digamma function or for the Hurwitz zeta-function (see  
Sections \ref{psi-tail} and \ref{Hurwitz-tail}).

Luckily, the other parts of the framework  described in the previous sections can be used for $G(z)$ and $\eta(s,z)$ too.
We will show the existence of suitable power series 
converging for every $z\in(0,2)$ and  using \eqref{difference-G-eta}
we will obtain useful results
on their reflection formulae as we did for the functions belonging to $\scrF$.

Following the previously mentioned ideas, in the next two subsections we will provide a detailed treatment 
for both $G(x)$ and  $\eta(s,x)$, $x>0$, $s>1$.

\subsection {The Bateman $G$-function}

The Bateman $G$-function is defined  as in \eqref{Bateman-G-def}
and  useful special values are
$G(1) 
= 2 \log 2$ and
$G(1/2) 
= \pi$.
According to Beebe \cite[p.~555]{Beebe2017},
it is not a good idea to compute $G(z)$ 
with its definition, 
and so we use the strategy described in Section \ref{intro-beyond-F}.
Our first result is the following
\begin{Proposition}[The number of summands for $G(x)$]
Let $x\in (0,1)$, $x\ne 1/2$,
$\acc\in \N$, $\acc\ge 2$,
\begin{equation}
\label{r-batemanG-estim}
r_G(x,\acc) = 
\max
\Bigl\{
\Bigl\lceil \frac{(\acc+1) \log 2 
+ 
\vert \log (1-\vert 1-x\vert)\vert + 1.2}{ \vert \log  \vert 1-x \vert \vert}\Bigr\rceil -1;
2
\Bigr\},
\end{equation}
and  $r^\prime_G(x,\acc) = r_G(1+x,\acc)$.
For $x\in (1/2,1)$ there exists $\theta_1=\theta_1(x)\in (-1/2,1/2)$ such that
\begin{equation}
\label{G-comp-x>1/2} 
G(x)   
= 2 \log 2 + \frac{2(1-x)}{x(1+x)} +  2 \sum_{k=1}^{r_G(x,\acc)} (1-2^{-k}) (\zeta(k+1)-1) (1-x)^{k}
+\vert \theta_1 \vert 2^{-\acc}. 
\end{equation}
For $x\in (0,1/2)$ we have that
there exists $\theta_2=\theta_2(x)\in (-1/2,1/2)$ such that
\begin{equation}  
\label{G-comp-x<1/2}
G(x)   
=  -2 \log 2 + \frac{2}{x} + \frac{2x}{(2+x)(1+x)} -  2 \sum_{k=1}^{r^\prime_G(x,\acc)} (1-2^{-k}) (\zeta(k+1)-1) (-x)^{k}
+\vert \theta_2 \vert 2^{-\acc}.   
\end{equation}
\end{Proposition}
\begin{Proof}
The starting point is that, using eq.~(1.17.6) on page 46 of 
\cite{ErdelyMOT1953}, we can write
\begin{equation}
\label{Bateman-series-def}
G(z)  =2 \log 2 +2\sum_{k=1}^{\infty} (1-2^{-k}) \zeta(k+1) (1-z)^{k},
\end{equation}
for every $z\in(0,2)$.
An alternative way to obtain \eqref{Bateman-series-def}
is inserting the duplication formula 
$2\psi(2w)= 2 \log 2 + \psi(w) + \psi(w+1/2)$
with $w=z/2$ into \eqref{Bateman-G-def}
thus getting 
\[
G(z) 
= 2 \log 2 + 2  \psi\Bigl(\frac{z+1}{2}\Bigr) - 2\psi(z)
\]
and remarking that the series of the right hand side is, for $z\in(0,2)$, the one in  \eqref{Bateman-series-def}.
Unfortunately, $(1-2^{-k}) \zeta(k+1)$ is not 
a decreasing sequence. However, subtracting  from \eqref{Bateman-series-def}
the difference of the geometric series of ratios $(1-z)$ and $(1-z)/2$, whose sum is 
 $\frac{1-z}{z} - \frac{1-z}{1+z} = \frac{1-z}{z(1+z)}$, it is easy 
to obtain, for $z\in(0,2)$, that
\begin{equation}
\label{Bateman-series-alt-def}
G(z)  =2 \log 2 + \frac{2(1-z)}{z(1+z)} + 2\sum_{k=1}^{\infty} (1-2^{-k}) (\zeta(k+1)-1) (1-z)^{k}.
\end{equation}
Letting $c_G(0):=2 \log 2 + \frac{2(1-z)}{z(1+z)} $ and
\begin{equation}
\label{Bateman-coeffs-def}
c_G(k) := 2(1-2^{-k}) (\zeta(k+1)-1)
\end{equation}  
for $k\ge 1$, and $C_G(k):=c_G(k)$, by Lemma \ref{elementary-estim} one gets
\[
\frac{c_G(k+1)}{c_G(k)}
=\frac{\zeta(k+2)-1}{\zeta(k+1)-1}\frac{1-2^{-k-1}}{1-2^{-k}}
<
\frac{1}{2}\frac{k+3}{k+1}  \frac{2^{k+1}-1}{2^{k+1}-2}
<1
\]
for $k\ge 2$, and this implies that point \ref{control-number-summands-f}) of Definition \ref{set-def} holds 
with $k_G=2$. 
Moreover, since the series in \eqref{Bateman-series-alt-def} absolutely converges for 
$z\in(0,2)$, we have that point \ref{series-f}) of Definition \ref{set-def} holds.
Since $c_G(k)$ is a decreasing sequence  
and $\vert \log(c_G(2))\vert = \vert  \log((3/2)(\zeta(3)- 1)) \vert < 1.2$, 
we obtain that that $r_G(z,\acc)$ can be estimated as in 
\eqref{r-batemanG-estim}.
Hence we can compute $G(x)$, $x\in (0,1)$, in the following way.   
For $x\in (1/2,1)$, using directly \eqref{Bateman-series-alt-def} and \eqref{r-batemanG-estim}, 
we  obtain \eqref{G-comp-x>1/2}.
For $x\in (0,1/2)$,
we first use \eqref{difference-G-eta} to write $G(x)= 2/x - G(1+x)$
and then we evaluate  \eqref{G-comp-x>1/2}  into $1+x$.
Defining $r^\prime_G(x,\acc) = r_G(1+x,\acc)$, formula \eqref{G-comp-x<1/2}
follows.
\end{Proof}

The reflection formulae for $G(x)$ follow from \eqref{G-comp-x>1/2}-\eqref{G-comp-x<1/2} and
\eqref{difference-G-eta} arguing as in the proof of Proposition \ref{DIF-formulae-f}.
\begin{Proposition}[The reflection formulae for $G(x)$]
Let $x\in (0,1)$, $x\ne 1/2$,
$\acc\in \N$, $\acc\ge 2$, $r_1(x,\acc)=r^\prime_G(x,\acc)/2$
and $r_2(x,\acc)  = r_G(x,\acc)/2$.
Recalling \eqref{Bateman-coeffs-def},
there exists $\theta=\theta(x)\in (-1/2,1/2)$  such that
for $0<x <1/2$ we have
\begin{align*}
G(x)+G(1-x) 
&=
\frac{2}{x}  +  \frac{2x}{(2+x)(1+x)} + \frac{2x}{(2-x)(1-x)} 
+  
2\sum_{\ell=1}^{r_1} c_G(2\ell-1) x^{2\ell-1}  
+ \vert \theta \vert 2^{-\acc},
\\
G(x)-G(1-x) 
&=
 - 4\log 2 + \frac{2}{x}  +  \frac{2x}{(2+x)(1+x)} - \frac{2x}{(2-x)(1-x)}
- 
2\sum_{\ell=1}^{r_1} c_G(2\ell) x^{2\ell} 
+ \vert \theta \vert 2^{-\acc} ,
\end{align*}
and for $1/2<x <1$ we have
\begin{align*}
G(x)+G(1-x) 
&=
\frac{2}{1-x}  +   \frac{2(1-x)}{(3-x)(2-x)} +  \frac{2(1-x)}{x(1+x)} 
+  
2\sum_{\ell=1}^{r_2} c_G(2\ell-1)  (1-x)^{2\ell-1} 
+ \vert \theta \vert 2^{-\acc} ,
\\ 
G(x)-G(1-x) 
&=
 4\log 2 - \frac{2}{1-x}  -   \frac{2(1-x)}{(3-x)(2-x)} + \frac{2(1-x)}{x(1+x)} 
+ 
2\sum_{\ell=1}^{r_2} c_G(2\ell)  (1-x)^{2\ell} 
+ \vert \theta \vert 2^{-\acc} .
\end{align*}
\end{Proposition}
Similar formulae hold for the infinite series too.

\subsubsection{The tail of the $G$-function}
Due to \eqref{difference-G-eta} we cannot directly use the general recursive formula \eqref{main-recursive-f}. 
But we can use twice the one for $\psi$ in \eqref{difference-psi} thus obtaining the following result.
\begin{Proposition}[Recursive formulae for $G(x)$]
\label{Bateman-recursive}
Let $\ell \in \N$. We have that 
\begin{equation}
\label{G-integral}
G(2\ell)  
= 
2  - 2 \log 2 + \sum_{j=1}^{\ell-1} \Bigl( \frac{1}{j+1/2}  - \frac{1}{j}\Bigr)
\quad
\textrm{and}
\quad
G(2\ell+1)  
=
2 \log 2 - 2 +\frac{1}{\ell} -  \sum_{j=1}^{\ell-1} \Bigl( \frac{1}{j+1/2}  - \frac{1}{j}\Bigr).
\end{equation}
Moreover,  for $x\in (2\ell, 2\ell+1)$ we have
\begin{equation}
\label{G-real-case1}
G(x) 
=
G\Bigl(2\Bigl\{\frac{x}{2}\Bigr\}  \Bigr)  
+
\sum_{j=0}^{\ell -1} 
\Bigl( 
\frac{1}{\{x/2\} + j+ 1/2} - \frac{1}{\{x/2\} + j} 
\Bigr)
\end{equation}
and, for $x\in (2\ell+1, 2\ell+2)$ we obtain
\begin{equation}
\label{G-real-case2}
G(x)
=
- G\Bigl(2\Bigl\{\frac{x}{2}\Bigr\} - 1\Bigr)  
+\frac{2}{2\{x/2\}-1} 
+
\sum_{j=0}^{\ell -1} 
\Bigl( 
\frac{1}{\{x/2\} + j+ 1/2} - \frac{1}{\{x/2\} + j} 
\Bigr).
\end{equation}
\end{Proposition}

\begin{Proof}
Using  \eqref{difference-psi}, we  obtain
\[
G(2\ell)  
=  
\psi\Bigl(\frac{3}{2}\Bigr) - \psi(1)  + \sum_{j=1}^{\ell-1} \Bigl( \frac{1}{j+1/2}  - \frac{1}{j}\Bigr) 
\]
and the first part of \eqref{G-integral} follows by recalling $\psi(1)=-\gamma$ and $\psi(3/2)= 2+ \psi(1/2) = 2 -2\log 2 - \gamma$.
The second part of \eqref{G-integral} is then obtained using its first part 
and $G(2\ell+1)  =  1/\ell - G(2\ell) $.
Let now $x\in (2\ell, 2\ell+1)$. 
We have $\lfloor x/2\rfloor  = \ell$,
$\{x/2\}  \in (0,1/2)$, $\{x/2 +1/2\} = \{x/2\}  +1/2  \in (1/2,1)$,
so that, using  \eqref{difference-psi}, we obtain
\[
G(x) = \psi \Bigl(\Bigl\{\frac{x}{2}\Bigr\} + \ell+ \frac{1}{2} \Bigr) 
- \psi \Bigl(\Bigl\{\frac{x}{2}\Bigr\} + \ell \Bigr) 
= 
\psi\Bigl(\Bigl\{\frac{x}{2}\Bigr\} + \frac{1}{2} \Bigr)  
- \psi\Bigl(\Bigl\{\frac{x}{2}\Bigr\}  \Bigr) 
+ 
\sum_{j=0}^{\ell -1} 
\Bigl( 
\frac{1}{\{x/2\} + j+ 1/2} - \frac{1}{\{x/2\} + j} 
\Bigr)  
\]
and \eqref{G-real-case1}  follows using the
definition of $G(x)$ in \eqref{Bateman-G-def}.
Let now $x\in (2\ell +1, 2\ell+2)$. We get
$\lfloor x/2\rfloor  = \ell$, $\lfloor (x+1)/2\rfloor  = \ell+1$,
$\{x/2\}  \in (1/2,1)$, $\{x/2 +1/2\} = \{x/2\}  -1/2  \in (0,1/2)$,
so that
\begin{align*}
\notag
G(x) &=  
\psi\Bigl(\Bigl\{\frac{x}{2}\Bigr\}   
+ \ell + \frac12 \Bigr) 
- \psi\Bigl(\Bigl\{\frac{x}{2}\Bigr\}
+ \ell\Bigr) 
= 
\psi\Bigl(\Bigl\{\frac{x}{2}\Bigr\}
+ \frac12\Bigr)  
- \psi\Bigl(\Bigl\{\frac{x}{2}\Bigr\} 
\Bigr)
+
\sum_{j=0}^{\ell -1} 
\Bigl( 
\frac{1}{\{x/2\} + j+ 1/2} - \frac{1}{\{x/2\} + j} 
\Bigr) 
\\& 
=
\psi\Bigl(\Bigl\{\frac{x}{2}\Bigr\}
- \frac12\Bigr)  
- \psi\Bigl(\Bigl\{\frac{x}{2}\Bigr\} 
\Bigr)
+\frac{1}{\{x/2\}-1/2} 
+
\sum_{j=0}^{\ell -1} 
\Bigl( 
\frac{1}{\{x/2\} + j+ 1/2} - \frac{1}{\{x/2\} + j} 
\Bigr)  
\end{align*}
and \eqref{G-real-case2} follows using the
definition of $G(x)$ in \eqref{Bateman-G-def}.
\end{Proof}

We remark that  for $x\in (2\ell, 2\ell+1)$ we have $2\{x/2\} \in (0,1)$ and  
that for $x\in (2\ell+1, 2\ell+2)$ we have  $2\{x/2\} -1 \in (0,1)$; hence, in both cases,
we can use \eqref{G-comp-x>1/2}-\eqref{G-comp-x<1/2}
to compute, respectively, $G(2\{x/2\})$ and $G(2\{x/2\} -1)$.  

Moreover, it is now clear that to evaluate the tails in \eqref{G-integral}-\eqref{G-real-case2} we can use  
twice the estimates in  \eqref{digamma-EM};  
the final error terms  for $G(x)$, $x>0$, will be bounded by  
$2E^{\psi}_{m}(\{x/2\})$, where  $E^{\psi}_{m}(v)$ is defined in 
\eqref{error-EM-psi}. Hence the value of $m_G$ has to be suitably modified
with respect to the one of $m_\psi$,
to ensure that the error term for $G$ will be smaller than the desired 
accuracy. If this is not possible, 
we can use the formulae in this paragraph to 
build a suitable horizontal shift as we did for $\psi$ in Section \ref{psi-tail}.

\subsection{The alternating Hurwitz zeta-function $\eta(s,x)$, $s > 1$ being fixed, $x>0$.}
The  alternating Hurwitz zeta-function is defined  as
in \eqref{alt-Hurwitz-def}.
We remark that 
$\eta (s,1) = (1-2^{1-s}) \zeta(s)$
and
$\eta (s,1/2) = 2^{s} \beta(s)$, for every $s>1$,
where $\beta(s)$ is the Dirichlet $\beta$-function defined in Section \ref{Dirichlet-beta-def-sec}.
Using \eqref{difference-G-eta} we also get $\eta(s,2) = 1-\eta(s,1) = 2^{1-s}\zeta(s)$
for $s>1$.
Unfortunately $\eta(s,\cdot)\not \in \scrF$, since
its functional equation is given by \eqref{difference-G-eta}
and hence point \ref{recursive-f}) of Definition \ref{set-def} does not hold.
However, in this case too we can use the strategy  described in Section \ref{intro-beyond-F}.
Our first result is the following
\begin{Proposition}[The number of summands for $\eta(s,x)$]
Let $s>1$ be fixed, $x\in (0,1)$, $x\ne 1/2$,
$\acc\in \N$, $\acc\ge 2$.
Let 
$c_\eta(s,0,z):= z^{-s}-2(1+z)^{-s}+(1-2^{1-s}) (\zeta(s)-1)$
and 
 $c_\eta(s,k):=  (1-2^{1-k-s})  (\zeta(k+s)-1)/(k B(s,k))$
for $k\ge 1$ and  $s > 1$, where \(
B(u,v) 
\)
is the Euler beta-function and $\zeta(u)$ is the Riemann zeta-function.
Let further
\begin{equation}
\label{r-alt-Hurwitz-estim}
r_{\eta}(s,x,\acc) = 
\max
\Bigl\{
\Bigl\lceil \frac{(\acc+1) \log 2 + \vert \log ( 1-\vert 1-x\vert ) \vert 
+ \vert \log c_{\eta}(s,\lceil s\rceil+1)\vert}
{ \vert \log  \vert 1-x \vert \vert}\Bigr\rceil -1;
\lceil s\rceil +1
\Bigr\},
\end{equation}
and 
$r^\prime_{\eta}(s,x,\acc) = r_{\eta}(s,1+x,\acc)$.
For $x\in (1/2,1)$ there exists $\theta_1=\theta_1 (s,x)\in (-1/2,1/2)$ such that
\begin{equation} 
\label{alt-zeta-Hurwitz-x>1/2}
\eta (s,x)  
=  
\frac{1}{x^s} -\frac{2}{(1+x)^s} 
+
\sum_{k=0}^{r_{\eta}(s,x,\acc)}  
c_\eta(s,k) (1-x)^{k}
+\vert \theta_1 \vert 2^{-\acc}.
\end{equation}
For $x\in (0,1/2)$ we have that
there exists $\theta_2=\theta_2(s,x)\in (-1/2,1/2)$ such that
\begin{equation} 
\label{alt-zeta-Hurwitz-x<1/2}
\eta (s,x)  
= 
\frac{1}{x^{s}}   - \frac{1}{(1+x)^s} +\frac{2}{(2+x)^s} 
-  \sum_{k=0}^{r^\prime_{\eta}(s,x,\acc)} 
c_\eta(s,k) (-x)^{k}+\vert \theta_2 \vert 2^{-\acc}.
\end{equation}
\end{Proposition}

\begin{Proof}
Let $s>1$ be fixed.
Using, for $w\in (-1,1)$ and $n\ge 1$, the Taylor series 
\[
(n+w)^{-s} = \sum_{k=0}^{\infty}  \frac{\Gamma(k+s)}{(k!) \Gamma(s)} (-w)^k n^{-s-k},
\]
and recalling that $\eta (s,1) = (1-2^{1-s}) \zeta(s)$, 
we obtain
\begin{align*}
\eta(s,1+w) 
&=
\sum_{n=1}^{\infty}  \frac{(-1)^{n-1}}{(n+w)^s}
=
\sum_{k=0}^{\infty}  \frac{\Gamma(k+s)}{(k!) \Gamma(s)} (-w)^k 
\eta(s+k,1)
=
\sum_{k=0}^{\infty} \frac{\Gamma(k+s)}{(k!) \Gamma(s)}  (1-2^{1-k-s}) \zeta(k+s) (-w)^{k}
\end{align*}
which, letting $z=1+w\in (0,2)$, becomes
\begin{equation}
\label{alt-Hurwitz-series-def}
\eta(s,z)  =  \sum_{k=0}^{\infty} \frac{\Gamma(k+s)}{(k!) \Gamma(s)}  (1-2^{1-k-s}) \zeta(k+s) (1-z)^{k}.
\end{equation}
An alternative way to obtain \eqref{alt-Hurwitz-series-def} is inserting the duplication formula 
$ \zeta(s,w) + \zeta (s, w+1/2) = 2^s \zeta(s,2w)$ with $w=z/2$ into \eqref{alt-Hurwitz-def}
thus getting 
\[
\eta(s,z) 
= \zeta(s,z) - 2^{1-s}  \zeta\Bigl(s,\frac{1+z}{2}\Bigr),
\]
and remarking that the series of the right hand side is, for $z\in(0,2)$, the one in  \eqref{alt-Hurwitz-series-def}.
Isolating in \eqref{alt-Hurwitz-series-def} the contribution of the series of  $z^{-s}-2(1+z)^{-s}$, we have
\begin{equation}
\label{alt-Hurwitz-series-def-conv}
\eta(s,z)  = \frac{1}{z^{s}}-\frac{2}{(1+z)^{s}}
+ \sum_{k=0}^{\infty} \frac{\Gamma(k+s)}{(k!) \Gamma(s)}  (1-2^{1-k-s}) (\zeta(k+s)-1) (1-z)^{k},
\end{equation}
for $s>1$ fixed, and $z\in(0,2)$.
Let 
$c_\eta(s,0,z):= z^{-s}-2(1+z)^{-s}+(1-2^{1-s}) (\zeta(s)-1)$ 
and $c_\eta(s,k):=  (1-2^{1-k-s})  (\zeta(k+s)-1)/(k B(s,k))$
for $k\ge 1$ and  $s > 1$, where \(
B(u,v) 
\)
is the Euler beta-function.
Letting further $C_\eta(s,k):=c_\eta(s,k)$,
by Lemma \ref{elementary-estim} one gets
\[
\frac{c_\eta(s,k+1)}{c_\eta(s,k)}
=\frac{\zeta(k+s+1)-1}{\zeta(k+s)-1} \frac{k B(s,k)}{(k+1) B(s,k+1)} \frac{1-2^{-k-s}}{1-2^{1-k-s}}
<
\frac{1}{2}
\frac{k+s+2}{k+s}
\frac{k+s}{k+1}
\frac{2^{k+s}-1}{2^{k+s}-2}  
<1
\]
for $k\ge s+1$, and this implies that point \ref{control-number-summands-f}) of Definition \ref{set-def} holds 
with $k_\eta(s)=\lceil s \rceil + 1$. 
Moreover, since the series in \eqref{alt-Hurwitz-series-def-conv} absolutely converges for 
$z\in(0,2)$, 
we have that point \ref{series-f}) of Definition \ref{set-def} holds.
Since $c_\eta(s,k)$ is a decreasing sequence  for $k\ge s+1$, 
we obtain that that $r_{\eta}(s,z,\acc)$ can be estimated as in 
\eqref{r-alt-Hurwitz-estim}.

Hence we can compute $\eta(s,x)$, $s>1$ being fixed, $x\in (0,1)$, in the following way.   
For $x\in (1/2,1)$  using directly \eqref{alt-Hurwitz-series-def-conv} and \eqref{r-alt-Hurwitz-estim}, 
we  obtain \eqref{alt-zeta-Hurwitz-x>1/2}.
For $x\in (0,1/2)$,
we first use \eqref{difference-G-eta} to write $\eta(s,x)= x^{-s} - \eta(s,1+x)$
and then we evaluate  \eqref{alt-zeta-Hurwitz-x>1/2}  into $1+x$.
Defining $r^\prime_{\eta}(s,z,\acc) = r_{\eta}(s,1+z,\acc)$, formula \eqref{alt-zeta-Hurwitz-x<1/2}
 follows.
\end{Proof}

The reflection formulae for $\eta (s,x)  \pm \eta (s,1-x)$, $s>1$ being fixed, $x\in (0,1)$, can be easily 
obtained  using \eqref{alt-zeta-Hurwitz-x>1/2}-\eqref{alt-zeta-Hurwitz-x<1/2} as we did in the previous sections
for $\log \Gamma(x)$, $\psi(x)$, $\zeta(s,x)$, $\zeta^\prime(s,x)$ and $G(x)$.

\subsubsection{The tail of the $\eta(s,\cdot)$-function, $s > 1$ being fixed}
Due to \eqref{difference-G-eta} we cannot directly use the general recursive formula \eqref{main-recursive-f}. 
But we can use twice the one for $\zeta(s,x)$ in  \eqref{difference-Hurwitz} thus obtaining the following result.
\begin{Proposition}[Recursive formulae for $\eta(s,x)$]
\label{eta-recursive}
Let $\ell \in \N$ and $s>1$. We have that 
\begin{equation}
\label{eta-integral-even}
\eta(s,2\ell)  
= 
\frac{\zeta(s)}{2^{s-1}}
+
\sum_{j=1}^{\ell-1} \Bigl( \frac{1}{(2j+1)^s}  - \frac{1}{(2j)^s}\Bigr)
\end{equation}
and
\begin{equation}
\label{eta-integral-odd}
\eta(s,2\ell+1)  
=
- 
\frac{\zeta(s)}{2^{s-1}}
+\frac{1}{(2 \ell)^s} 
-
\sum_{j=1}^{\ell-1} \Bigl( \frac{1}{(2j+1)^s}  - \frac{1}{(2j)^s}\Bigr).
\end{equation}
Moreover,  for $x\in (2\ell, 2\ell+1)$ we have
\begin{equation}
\label{eta-real-case1}
\eta(s,x) 
=
\eta\Bigl(s,2\Bigl\{\frac{x}{2}\Bigr\}  \Bigr)  
+
\sum_{j=0}^{\ell -1} 
\Bigl( 
\frac{1}{(2\{x/2\} + 2j+ 1)^s} - \frac{1}{(2\{x/2\} + 2j)^s} 
\Bigr) 
\end{equation}
and, for $x\in (2\ell+1, 2\ell+2)$ we obtain
\begin{equation}
\label{eta-real-case2}
\eta(s,x)
=
- \eta\Bigl(s,2\Bigl\{\frac{x}{2}\Bigr\} - 1\Bigr)  
+\frac{1}{(2\{x/2\}-1)^s} 
+
\sum_{j=0}^{\ell -1} 
\Bigl( 
\frac{1}{(2\{x/2\} + 2j+ 1)^s} - \frac{1}{(2\{x/2\} + 2j)^s} 
\Bigr).
\end{equation}
\end{Proposition}

\begin{Proof}
Using  \eqref{difference-Hurwitz}, we obtain
\[
\eta(s,2\ell)  
=  
\frac{\zeta(s,1) - \zeta\bigl(s, 3/2 \bigr)}{2^s}  
+
\sum_{j=1}^{\ell-1} \Bigl( \frac{1}{(2j+1)^s}  - \frac{1}{(2j)^s}\Bigr) 
\]
and \eqref{eta-integral-even} follows by recalling 
\eqref{alt-Hurwitz-def} and  $\eta(s,2) = 2^{1-s}\zeta(s)$.
Formula \eqref{eta-integral-odd} is then obtained using   \eqref{eta-integral-even}
and $\eta(s,2\ell+1)  =  (2\ell)^{-s} - \eta(s,2\ell) $.
Let now $x\in (2\ell, 2\ell+1)$. 
We have $\lfloor x/2\rfloor  = \ell$,
$\{x/2\}  \in (0,1/2)$, $\{x/2 +1/2\} = \{x/2\}  +1/2  \in (1/2,1)$,
so that, using  \eqref{difference-Hurwitz}, we obtain
\begin{align}
\notag
\eta(s,x) &= 
2^{-s}\zeta\Bigl(s,\Bigl\{\frac{x}{2}\Bigr\} + \ell \Bigr)  
-
2^{-s} \zeta\Bigl(s,\Bigl\{\frac{x}{2}\Bigr\} + \ell+ \frac{1}{2} \Bigr) 
\\&
\notag
= \frac{\zeta(s, \{x/2\}) - \zeta(s, \{x/2\}+1/2)}{2^s}
+\sum_{j=0}^{\ell -1} 
\Bigl( 
\frac{1}{(2\{x/2\} + 2j+ 1)^s} - \frac{1}{(2\{x/2\} + 2j)^s} 
\Bigr) 
\end{align}
and \eqref{eta-real-case1}  follows using the
definition of $\eta(s,x)$ in \eqref{alt-Hurwitz-def}.
Let now $x\in (2\ell +1, 2\ell+2)$. We get
$\lfloor x/2\rfloor  = \ell$, $\lfloor (x+1)/2\rfloor  = \ell+1$,
$\{x/2\}  \in (1/2,1)$, $\{x/2 +1/2\} = \{x/2\}  -1/2  \in (0,1/2)$,
so that
\begin{align}
\notag
\eta(s,x)  
&=
2^{-s}\zeta\Bigl(s,\Bigl\{\frac{x}{2}\Bigr\} + \ell \Bigr)  
-
2^{-s} \zeta\Bigl(s,\Bigl\{\frac{x}{2}\Bigr\} + \ell+ \frac{1}{2} \Bigr)  
\\&
\notag
=
\eta\Bigl(s,2\Bigl\{\frac{x}{2}\Bigr\} \Bigr)  
+
\sum_{j=0}^{\ell -1} 
\Bigl( 
\frac{1}{(2\{x/2\} + 2j+ 1)^s} - \frac{1}{(2\{x/2\} + 2j)^s} 
\Bigr)  
\\&
\notag
=
- \eta\Bigl(s,2\Bigl\{\frac{x}{2}\Bigr\} - 1\Bigr)  
+\frac{1}{(2\{x/2\}-1)^s} 
+
\sum_{j=0}^{\ell -1} 
\Bigl( 
\frac{1}{(2\{x/2\} + 2j+ 1)^s} - \frac{1}{(2\{x/2\} + 2j)^s} 
\Bigr)  ,
\end{align}
and \eqref{eta-real-case2}  follows using the
definition of $\eta(s,x)$ in \eqref{alt-Hurwitz-def}.
\end{Proof}

We remark that  for $x\in (2\ell, 2\ell+1)$ we have $2\{x/2\} \in (0,1)$ and  
that for $x\in (2\ell+1, 2\ell+2)$ we have  $2\{x/2\} -1 \in (0,1)$; hence, in both cases,
we can use \eqref{alt-zeta-Hurwitz-x>1/2}-\eqref{alt-zeta-Hurwitz-x<1/2}
to compute, respectively, $\eta(s,2\{x/2\})$ and $\eta(s,2\{x/2\} -1)$.

Moreover, it is now clear that to evaluate the tails in \eqref{eta-integral-even}-\eqref{eta-real-case2} we can 
use twice equation  \eqref{Hurwitz-EM};  
the final error terms  for $\eta(s,x)$, $x>0$, $s>1$, will be bounded by  
$2^{1-s}E^{\zeta_H}_{m}(s,\{x/2\})$, where  $E^{\zeta_H}_{m}(s,v)$ is defined in 
\eqref{error-EM-Hurwitz}. Hence the value of $m_\eta$ has to be suitably modified
with respect to the one of $m_{\zeta_H}$,
to ensure that the error term for $\eta(s,x)$ will be smaller than the desired 
accuracy. If this is not possible, 
we can use the formulae in this paragraph to 
build a suitable horizontal shift as we did for $\zeta(s,x)$ in Section \ref{Hurwitz-tail}. 

\section{Some practical tests using Pari/GP}
\label{tests}

We performed some tests by implementing in Pari/GP \cite{PARI2021} the algorithms here described\footnote{
The programs can be downloaded at the  
page: \url{http://www.math.unipd.it/~languasc/specialfunctions.html};
some practical examples are included towards the bottom of the scripts.}.
The running time comparisons refer to the internal Pari/GP
functions (however, such functions work for complex inputs too
while ours work only for positive inputs) and are
obtained on a Dell OptiPlex-3050, 
equipped with an Intel i5-7500 CPU, 3.40GHz, 16GB of RAM and running 
Ubuntu 20.04.3LTS. 
A professional implementation of our algorithm should be able to obtain better practical performances
than the ones described below.

We decided to use Pari/GP to perform such tests since it allows easily to 
work with a multiprecision library and many of the functions here mentioned
are already implemented.
Clearly, many other implementations are available; our main goal
here is to show that our unified computational strategy can represent 
a possible alternative and that, in some particular cases, like the 
ones involving the reflection formulae, it
has a remarkable efficiency.

We will allow an accuracy of $32,64,80$ or $128$ bits in our implementations
of these functions but it is easy to modify this to allow an arbitrarily large precision.
The choice of the length of the sum, i.e., of the parameter $m$, in the Euler-Maclaurin formula is performed at 
runtime, even if for  $\log\Gamma(x)$, $\psi(x)$ and $G(x)$ it is also possible to use Tables \ref{Gamma-tail-accuracy-table}
and \ref{psi-tail-accuracy-table}.
The precomputation of the needed Riemann zeta-function  values
gives us excellent performances;
in particular, for $\log\Gamma(x),\psi(x), G(x)$ this works nicely for $x\in(0,1)$ and in this case
our algorithm is about twice times faster that the internal Pari/GP functions.\footnote{On 
November 8th, 2021, K.~Belabas \cite{Belabas2021},
maintainer and developer of Pari/GP, communicated me that some of the
ideas used here and in \cite{Languasco2021} to compute $\Gamma(x)$ and
$\log\Gamma(x)$ for $x\in(0,1)$ have been used to improve
their Pari/GP implementations (from version 2.9.0 on) and that further tunings of their code
are about to be released. Analogous ideas
for $\psi(x)$ will be inserted in Pari/GP-2.15.
(Added on November 13th, 2022: in fact such tunings were, at least for the $\log\Gamma$-function, 
included in Pari/GP-2.15).}
For $x>1$ the computation of the tail of such functions requires too much time and 
the performances of our method become worse than the ones used in Pari/GP;
in this case it might be better to use an asymptotic formula
instead of the recursive step of point \ref{recursive-f})
of Definition \ref{set-def}.

We also tested the case of the Hurwitz zeta-function by computing $10000$ times
the same values with our algorithm and the internal Pari/GP functions (both
with an accuracy of $128$ bits). For integral values of $s$, our script is about $5.5$-times faster and
for non-integral values of $s$ the performances are even better; for example,
we evaluated  $\zeta(8.3,1345.1234)$ for $10000$ times in less than $4.9$ seconds with our 
algorithm while the internal functions of Pari/GP required  about two minutes
and $21$ seconds of time (with a factor of about $30$ as a performance gain). 
To test the implementation designed  for the Fast Fourier Transform applications, see Section \ref{DirichletL}, 
we have chosen a small prime $q=305741$ and we computed $\zeta(s,a_k/q)$ 
and $\zeta^{\prime}(s,a_k/q)$ for some values of $s$, where $a_k \equiv g^k \bmod{q}$, 
$g$ is a primitive root of $q$ and $k$ runs from $0$ to $q-2$. 
In the worst case our script is twice times faster than
the internal Pari/GP functions. 
Moreover, our algorithm seems to be particularly efficient
for  $\zeta^{\prime}(s,a_k/q)$: for example, with an accuracy of $128$ bits,
$s=8.3$,  $q=307541$ and $k$ running from $0$ to $q-2$, our implementation is about $32$-times 
faster than the one that uses directly the internal Pari/GP functions. 
Moreover, this can be further improved using the script that computes
at the same time both $\zeta(s,a_k/q)$ and $\zeta^{\prime}(s,a_k/q)$.
We also remark that
computing $\zeta(s,a_k/q)\pm \zeta(s,1-a_k/q)$ and $\zeta^{\prime}(s,a_k/q)\pm \zeta^{\prime}(s,1-a_k/q)$
gives much better performances thanks to the use of the reflection formulae and
Section \ref{further-FFT-gain}: with an accuracy of $128$ bits,
$s=8.3$,  $q=307541$ and $k$ running from $0$ to $q-2$, our implementation is about $60$-times 
faster than the one that uses directly the internal Pari/GP functions. 
For $q=6766811$ the gain factor in the speed of the last computation is about $72$-times;
in fact, we think that the larger $q$ we use, the better the performance gain  becomes.
It is hence clear that in computing the values of the Dirichlet $L$-functions, see Section \ref{DirichletL},
the use of our algorithm gives excellent performances in computing the input sequences
of the Fast Fourier Transform procedures, see Section \ref{DIF-sect} about the
decimation in frequency strategy.

For the Dirichlet $\beta$-function we computed $10000$ times
the same values both with our algorithm and the internal Pari/GP functions (both
with a precision of $128$ bits).
For integral values of $s$, our script is about three times faster and
for non-integral values of $s$ the performances are even better; for example 
we evaluated  $\beta(8.3)$ for $10000$ times in $535$ milliseconds with our 
algorithm while the internal functions of Pari/GP required about seven seconds of time.
A better performance gain is obtained for both $\beta^\prime(s)$ 
and $\beta^\prime(s)/\beta(s)$; with $s=8.3$, our algorithm is  
respectively faster by a factor of $39$ and $52$. Moreover, using our script for  $\beta^\prime(s)/\beta(s)$ 
is particularly convenient since is also gives the values of $\beta^\prime(s)$ and $\beta(s).$

From all the previous tests and examples, we can conclude that the algorithm here presented
is particularly useful when it is possible to exploit the precomputation
of the $c_f(s,k)$-coefficients of the series in point \ref{series-f}) of Definition \ref{set-def}, 
as for the computation of the values of $\zeta (s,a_k/q)$ and $\zeta^\prime (s,a_k/q)$ needed
to obtain $L(s,\chi)$ and $L^\prime(s,\chi)$, where $s>1$, $\chi$ 
runs over the non-principal Dirichlet character modulo $q$, for every $q$ in a large set of odd prime numbers.
\renewcommand{\bibliofont}{\normalsize}

\vskip 0.5cm
\noindent
Alessandro Languasco,
Universit\`a di Padova,
Dipartimento di Matematica,
``Tullio Levi-Civita'',
Via Trieste 63,
35121 Padova, Italy.    
{\it e-mail}: alessandro.languasco@unipd.it   

\end{document}